\documentclass[12pt]{amsart}

\usepackage{graphicx}
\usepackage{times}
\newtheorem{thm}{Theorem}[section]
\newtheorem{cor}[thm]{Corollary}
\newtheorem{notation}[thm]{Notation}

\newtheorem{defin}[thm]{Definition}
\newtheorem{lemma}[thm]{Lemma}

\newtheorem{prop}[thm]{Proposition}
\newtheorem{rmk}[thm]{Remark}

\newcommand{\aaa}{\mbox{$\alpha$}}

\newcommand{\mcA}{\mbox{$\mathcal A$}}
\newcommand{\mcB}{\mbox{$\mathcal B$}}

\newcommand{\bbb}{\mbox{$\beta$}}

\newcommand{\Sss}{\mbox{$\Sigma$}}

\newcommand{\Ggg}{\mbox{$\Gamma$}}
\newcommand{\Lll}{\mbox{$\Lambda$}}

\newcommand{\bdd}{\mbox{$\partial$}}

\newcommand{\xX}{\mbox {\sc x}}
\newcommand{\yY}{\mbox {\sc y}}
\newcommand{\aA}{\mbox {\sc a}}
\newcommand{\bB}{\mbox {\sc b}}

\begin{document}



\title{Multiple bridge surfaces restrict knot distance}

\author{Maggy Tomova}
\address{\hskip-\parindent
 Maggy Tomova\\
  Mathematics Department \\
University of Iowa \\
 Iowa city, IA 52240, USA}
\email{mtomova@math.uiowa.edu}

\thanks{Research partially supported by an NSF grant.}

\date{\today}

\begin{abstract}
    
  Suppose $M$ is a closed irreducible orientable $3$-manifold, $K$ is a 
  knot in $M$, $P$ and $Q$ are bridge surfaces for $K$ and $K$ is not removable 
 with respect to $Q$.  We
  show that either $Q$ is equivalent to $P$ or $d(K,P)\leq 2-\chi(Q-K)$.  
  If $K$ is not a 2-bridge knot, then the result
  holds even if $K$ is removable with respect to $Q$. As a corollary we show that if a
  knot in $S^3$ has high distance with respect to some bridge sphere
  and low bridge number, then the knot has a unique minimal bridge
  position.
    
\end{abstract}

\maketitle

\section{Introduction and Preliminaries} \label{sec:intro}
Distance is a generalization of the concepts of weak and strong
compressibility for bicompressible surfaces originally due to Hempel
\cite{Hem}.  It has been successfully applied to study Heegaard
splittings of 3-manifolds.  For example in \cite{Hart} Hartshorn shows
that the Euler characteristic of an essential surface in a manifold
bounds the distance of any of its Heegaard splittings.  In \cite{ST},
Scharlemann and Tomova show that the Euler characteristic of any
Heegaard splitting of a 3-manifold similarly bounds the distance of
any non-isotopic splitting.

A knot $K$ in a 3-manifold $M$ is said to be in bridge position with 
respect to a surface $P$ if $P$ is a Heegaard surface for $M$ and $K$ 
intersects each of the components of $M-P$ in arcs that are parallel 
to $P$. If $K$ is in bridge position 
with respect to $P$, we say that $P$ is a bridge surface for $K$.
The definition of distance has been extended to apply to bridge 
surfaces. In \cite{BS}, Bachman and Schleimer prove that Hartshorn's result
extends to the distance of a bridge surface, namely the Euler characteristic of
an essential properly embedded surface in the complement of a knot
bounds the distance of any bridge surface for the knot. 
In this paper we extend the ideas in \cite{ST} to show that the result
there also extends to the case of a knot with two different bridge surfaces.

\medskip

{\bf Theorem:} {\em Suppose $K$ is a non-trivial knot in a closed, irreducible and
orientable $3$-manifold $M$ and $P$ is a bridge surface for $K$ that 
is not a 4-times punctured sphere. If 
$Q$ is also a bridge surface for $K$ that is not equivalent to $P$, 
or  
if $Q$ is a Heegaard surface for $M-\eta (K)$ then
$d(K,P) \leq 2- \chi(Q-K)$.}

\medskip 

Two Heegaard splittings are usually considered to be equivalent if one 
is isotopic to a possibly stabilized copy of the other. For bridge 
surfaces there are three obvious geometric operations that correspond 
to stabilizations and they are described in Section \ref{sec:main}.

A knot $K$ is said to be removable with respect to a bridge surface 
$Q$ if $K$ can be isotoped to lie in the spine of one of the 
handlebodies $M-Q$. Thus after the isotopy, $Q$ is a Heegaard surface 
for $M-\eta(K)$. If we restrict our attention only to bridge surfaces with 
respect to which the knot is not removable, we may extend the above 
theorem also to 2-bridge knots.

\medskip

{\bf Corollary:} {\em
    Suppose $P$ and $Q$ are two bridge surfaces for a knot $K$ and $K$ 
    is not removable with respect to $Q$. Then either $Q$ is 
    equivalent to $P$ or $d(P)\leq 2-\chi(Q_K)$.} 

    \medskip

The result proves a conjecture of Bachman and Schleimer put forth in
\cite{BS}.

 \medskip
 
  {\bf Corollary:} {\em If $K\subset S^3$ is in minimal bridge
  position with respect to a sphere $P$ such that $d(K,P) > |P\cap K|$
  then $K$ has a unique minimal bridge position.}
 
 \medskip

 The basic idea of the proof of the above theorem is to consider a
 2-parameter sweep-out of $M-K$ by the two bridge
 surfaces. We keep track of information about compressions by introducing labels for the
 regions of the graphic associated to the sweep-out.  We are able to
 conclude that if particular combinations of labels occur we can
 deduce the desired result.  Using a quadrilateral version
 of Sperner's lemma, we conclude that one of the label
 combinations we have already considered must occur.

  \section{Surfaces in a handlebody intersected by the knot in
  unknotted arcs} \label{sec:surfaceinKhandlebody}
  
  Throughout this paper we will use the
  following definitions and notation:
  
  \begin{notation}
      Let $M$ be a compact orientable irreducible 3-manifold. 
      If $K \subset M$ is some properly embedded 1-manifold, let $M_K$
   denote $M$ with a regular (open) neighborhood $N(K)$ of $K$ removed.  If
   $X$ is any subset of $M$, let $X_K=M_K\cap X$.
      
  \end{notation}

\begin{defin} \label{def:long} Suppose $(F,\bdd F)\subset (M, \bdd M)$
is a properly embedded surface in a compact orientable irreducible manifold $M$ containing a
1-manifold $K$ such that $F$ is transverse to $K$.
\begin{itemize}

    \item We will say that $F_K$ is \underline{n-times punctured} if
    $|F\cap K|=n$.  If $F_K$ is 1-time punctured, we will call it
    \underline{punctured}.
    
\item A simple closed curve on $F_K$ is \underline{inessential} if it
bounds a subdisk of $F_K$ or it is parallel to a component of $\bdd
F_K$.  Otherwise the curve is \underline{essential}.

\item A properly embedded arc $(\bbb, \bdd \bbb) \subset(F_K,\bdd
F_K)$ is \underline{essential} if no component of $F_K-\bbb$ is a
disk.

\item A properly embedded disk $(D, \bdd D)\subset (M_K,F_K)$ is a
\underline{compressing disk} for $F_K$ in $M_K$ if $\bdd D$ is an
essential curve in $F_K$.

\item A \underline{cut disk} $D^c$ for $F_K$ is a punctured disk such
that $\bdd D^c =\alpha \cup \beta$ where $D^c \cap cl(N(K))=\alpha$, $D^c
\cap F_K=\beta$ and $\beta$ is an essential curve in $F_K$.

\item A \underline{c-disk} $D^*$ for $F_K$ is either a cut-disk or a
compressing disk.

\item A surface $F_K$ is called \underline{incompressible} if it has
no compressing disks, \underline{cut-incompressible} if it has no 
cut-disks and \underline{c-incompressible} if it has no c-disks.

\item A surface $F_K$ is called \underline{essential} if it is
incompressible and at least one of its components is not parallel to
$\bdd M_K$.
\end{itemize}

\end{defin}

Now we restrict our attention to the case when the 3-manifold we are
considering is a handlebody and the 1-manifold $K$ consists of
``unknotted" properly embedded arcs.  To make this more precise we use
the following definition modeled after the definition of a
$K$-compression body introduced in \cite{Bach}.

  \begin{defin}
    A $K$-handlebody, $(A,K)$ is a handlebody $A$ and a 1-manifold,
    $(K,\bdd K) \subset (A, \bdd A)$, such that $K$ is a disjoint
    union of properly embedded arcs and for each arc $\kappa \in K$
    there is a disk, $D \subset A$ with $\bdd D=\kappa \cup \alpha$,
    where $D \cap K=\kappa$ and $D \cap \bdd A = \alpha$.  These disks
    are called bridge disks.
     
  \end{defin}
  
  Many results about handlebodies have analogues for $K$-handlebodies. 
  We will need some of these.
  
  A {\em spine} $\Sigma_A$ of a handlebody $A$ is any graph that $A$
  retracts to.  Removing a neighborhood of a spine
  from a handlebody results in a manifold that is homeomorphic to
  $surface \times I$.  We need a similar
  concept for a spine of a $K$-handlebody.  Following \cite{BS} we
  define the spine $\Sigma_{(A,K)}$ of the $K$-handlebody $(A,K)$ to be the union of a
  spine of the handlebody $A$, $\Sigma_A$, together with a collection of
  straight arcs $t_i$, one for each component of $K \cap A$, where 
  one endpoint of each $t_i$ lies on $\Sigma_A$ and the other 
  endpoints of the $t_i$ lie on distinct components of $K \cap A$.  If
  $\bdd A=P$ then $A_K -\Sigma_{(A,K)}\cong P_K \times I$.  As in the
  handlebody case, spines of $K$-handlebodies are not unique.

\begin{notation}
    For the rest of this paper, unless otherwise specified, let
    $(A,K)$ be a $K$-handlebody with $P=\bdd A$ and spine
    $\Sigma_{(A,K)}$.  We will always assume that if $A$ is a ball, 
    then $K$ has at least 3 components. 
     $F \subset A$ will be a properly embedded surface that
    is transverse to $K$.  We continue to denote by $N(K)$ a regular
    neighborhood of $K$.
\end{notation}

\begin{defin}
Two embedded meridional surfaces $S$ and $T$ in $(M,K)$ are called
$K$-parallel if they cobound a region homeomorphic to $S_K \times I$
i.e. the region of parallelism contains only unknotted segments of $K$
each with one endpoint on $S$ and one endpoint on $T$.

Two meridional surfaces $S$ and $T$ are $K$-isotopic if there exists an isotopy
from $S$ to $T$ so that $S$ remains transverse to $K$ throughout the 
isotopy. 

\end {defin}

\begin{lemma} \label{lem:removable}
If $(E, \bdd E)\subset (A_K, P_K)$ is a possibly punctured disk such
that $\bdd E$ is an inessential curve on $P_K$, then $E$ is parallel
to a possibly punctured subdisk of $P_K$.
\end{lemma}

\begin{proof}

Let $E'$ be the possibly punctured disk $\bdd E$ bounds on $P_K$. 
There are three cases to consider.  If $E$ and $E'$ are both disks,
then they cobound a ball as $A_K$ is irreducible, and thus $E$ is
parallel to $E'$.  If one of $E$ and $E'$ is a once punctured disk and
the other one is a disk, then the sphere $E \cup E'$ intersects $K$ only once.  
The manifold is irreducible and $E \cup E'$ is
separating so this is not possible.  Finally, if both $E$ and $E'$ are
once punctured disks, then by irreducibility of $A$ and the definition
of a $K$-handlebody, $E$ and $E'$ cobound a product region in $A_K$. 
This product region intersects some bridge disk for $K$ in a single
arc, so the arc of $K$ between $E$ and $E'$ is a product arc.  It
follows that $E$ and $E'$ are parallel as punctured disks.

\end{proof}

 \begin{defin} \label{defin:boundcomp}
A $P$-compressing disk for $F_K\subset A_K$ is a disk $D \subset A_K$
so that $\bdd D$ is the end-point union of two arcs, $\aaa = D \cap
P_K$ and $\bbb = D \cap F_K$, and $\bbb$ is an essential arc in $F_K$.

 \end{defin}
 
 The operation of compressing, cut-compressing and $P$-compressing the surface
 $F_K$ have natural duals that we will refer to as {\em tubing} 
 (possibly tubing along a subset of the knot) and
 {\em tunneling} along an arc dual to the c-disk or the $P$-compressing disk.  The precise definitions of these operations
were given in \cite{Sc2}: Suppose $F\subset M$ is a properly embedded 
surface in a manifold containing a knot $K$. Let $\gamma\subset 
interior(M)$ be an embedded arc such 
that $\gamma \cap F = \bdd \gamma$. There is a relative tubular 
neighborhood $\eta (\gamma) \cong \gamma \times D^2$ so that $\eta 
(\gamma)$ intersects $F$ precisely in the two diskfibers at the ends 
of $\gamma$. Then the surface obtained from $F$ by removing these two 
disks and attaching the cylinder $\gamma \times \bdd D^2$ is said to 
be obtained by tubing along $\gamma$.  We allow for the possibility 
that $\gamma \subset K$.
Similarly if $\gamma \subset \bdd M$, there is a relative 
neighborhood $\eta (\gamma) \cong \gamma \times D^2$ so that $\eta 
(\gamma)$ intersects $F$ precisely in the two diskfibers at the ends 
of $\gamma$ and $\eta(\gamma)$ intersects $\bdd M$ in a
rectangle. Then the surface obtained from $F$ by removing the two 
half disks and attaching the rectangle $(\gamma \times \bdd D^2) \cap 
M$ is said to 
be obtained by tunnelling along $\gamma$.

    We will have many occasions to use $P$-compressions of surfaces so
    we note the following lemma.

  \begin{lemma} \label{lemma:boundredessential}
   Suppose $F_K \subset A_K$ is a properly embedded surface and $F'_K$
   is the result of $P$-compressing $F_K$ along a $P$-compressing disk
   $E_0$.  Then
    \begin{enumerate}
    \item If $F_K'$ has a c-disk, $F_K$ also has a c-disk of the same
    kind (cut or compressing).  
    \item If $F_K$ intersects every spine $\Sigma_{(A,K)}$ then so does $F_K'$.  
    \item Every curve of $\bdd F_K$ can be isotoped on $P_K$ to be disjoint from any curve
    in $\bdd F_K'$.

\end{enumerate}

   \end{lemma}

   \begin{proof}

       The original surface $F_K$ can be recovered from $F_K'$ by
       tunneling along an arc that is dual to the $P$-compressing
       disk.  This operation is performed in a small neighborhood of
       $P_K$ so if $F_K'$ has compressing or cut disks, they will be
       preserved in $F_K$.  Also if $F_K'$ is disjoint from some
       $\Sigma_{(A,K)}$, then adding a tunnel close to $P_K$ will not
       introduce any intersections with this spine.
       
       For the last item consider the frontier of $N(F_K\cup E_0) \cap
       P_K$ where $N$ denotes a regular neighborhood.  This set of
       disjoint embedded curves on $P_K$ contains both $F_K\cap P_K$
       and $F_K'\cap P_K$.

       \end{proof}

In the case of a handlebody it is also known that any essential
surface must have boundary.  The following lemma proves the
corresponding result for a $K$-handlebody.
 
     \begin{lemma} \label{lem:notclosed}
    If $F_K$ is an incompressible surface in $A_K$, then one of the 
    following holds,
    \begin{enumerate}
\item $F_K$ is a sphere,
\item $F_K$ is a twice punctured sphere, or
\item $F_K\cap P_K \neq \emptyset$.
\end{enumerate}
 \end{lemma}

      \begin{proof}
 Suppose $F_K$ is an incompressible surface in $A_K$ that is not a 
 sphere or a twice-punctured sphere, such that
 $P_K\cap F_K=\emptyset$.  Let $\Delta$ be the collection of a
 complete set of compressing disks for the handlebody $A$ together
 with all bridge disks for $K$.  Via an innermost disk argument,
 using the fact that $F_K$ is incompressible, we may assume that $F_K
 \cap \Delta$ contains only arcs.  Any arc of intersection between a
 disk $D \in \Delta$ and $F_K$ must have both of its endpoints lying
 on $N(K)$ as $F_K \cap P_K = \emptyset$ and thus lies on one of
 the bridge disks.  Consider an outermost such arc on $D$ cutting a
 subdisk $E$ of $D$.  Doubling $E$ along $K$ produces a compressing 
disk for
 $F_K$ which was assumed to be incompressible.  Thus $F_K$ must be
 disjoint from $\Delta$ and therefore $F_K$ lies in the ball $A_K -
 \Delta$ contradicting the incompressibility of $F_K$.
 \end{proof}

 Finally it is well known that if $F$ is a closed connected 
incompressible surface 
 contained in $A-\Sigma_A \cong P\times I$, then $F$ is isotopic to $P$. A 
similar 
 result holds if we consider $F_K\subset (A_K-\Sigma_{(A,K)})=P_K\times 
I$.

 \begin{lemma}\label{lem:collar} 
     Suppose $P$ is a closed connected surface, and $K\neq \emptyset$ is a 1-manifold 
    properly embedded in $P \times I$ so that each component of $K$ can be 
     isotoped to be vertical with respect to the product structure.
     If $F_K\subset P_K\times I$ is a properly embedded connected incompressible 
     surface such that $F_K \cap (P \times \{0\})=F_K \cap (P \times 
     \{1\})=\emptyset$, then one of the following holds,
     \begin{enumerate}
	 \item $F_K$ is a sphere disjoint from the knot,
\item $F_K$ is a twice punctured sphere, or
\item $F_K$ is $K$-isotopic to $P_K\times \{0\}$. 
\end{enumerate}	
 \end{lemma}    
    
      \begin{proof}
	  Suppose $F_K$ is not a sphere or a twice punctured sphere. 
	  Consider the set $S$ consisting of properly embedded arcs on 
	   $P_K$ so that $P_K-S$ is a disk. This collection gives rise to a 
	   collection $\Delta=S\times I$ of disks in $P_K\times I$ so 
	   that $(P_K\times I)-\Delta$ is a ball.  As $F_K$ is 
incompressible, 
	   by an innermost disk argument we may assume that it does not 
	   intersect $\Delta$ in any closed curves. If $F_K \cap \Delta$ 
	   contains an arc that has both of its endpoints on the same 
	   component of $K$, doubling the subdisk of $\Delta$ an outermost 
	   such arc bounds would give a compressing disk for $F_K$. 
	   Consider the components of $F_K$ lying in the ball 
$(P_K\times I)-\Delta$. As $F_K$ is incompressible all of these components must be 
	   disks. In fact, as $F_K$ is connected, there is a single disk 
	   component. This disk is isotopic to $(P_K-S)\times 0$ 
and the 
	   maps that glue $(P_K\times I)-\Delta$  to recover $P_K\times I$ do not 
	   affect the isotopy.

 \end{proof}

 \section{The curve complex and distance of a knot} 
 \label{sec:curvecomplex}

 Suppose $V$ is a compact, orientable, properly embedded surface in a 
3-manifold
 $M$.  The {\em curve complex} of $V$ is a graph $\mathcal{C}(V)$, 
with vertices corresponding to isotopy
 classes of essential simple closed curves on $V$.  Two vertices are
 adjacent if their corresponding isotopy classes of curves have
 disjoint representatives.  If $S$ and $T$ are subsets of vertices of
 $\mathcal{C}(V)$, then $d(S,T)$ is the length of the shortest path in
 the graph connecting a vertex in $S$ and a vertex in $T$.

 \begin{defin}
      Let $(P,\bdd P)\subset(M,\bdd M)$ be a properly embedded 
surface 
      in an orientable irreducible 3-manifold $M$. $P$ will be called 
      a splitting surface if $M$ is the union of two manifolds $A$ 
      and $B$ along $P$. We will say $P$ splits $M$ into $A$ and 
$B$.  
  \end{defin}

 If $P$ is a closed embedded bicompressible surface with 
 $\chi(P)<0$ splitting $M$ into
 submanifolds $A$ and $B$, let $\mcA$ (resp $\mcB$) be the set of all
 simple closed curves on $P$ that bound compressing disks for $P$ in
 $A$ (resp $B$).  Then $d(P)=d(\mcA,\mcB)$ i.e, the length of the
 shortest path in the graph $\mathcal{C}(P)$ between a curve in 
$\mcA$ and a curve in
 $\mcB$.  If $d(P) \leq 1$, i.e. there are compressing disks on
 opposite sides of $P$ with disjoint boundaries, then the surface $P$
 is called {\em strongly compressible} in $M$.  Otherwise $P$ is {\em
 weakly incompressible}.

 Much like bridge number and width, the distance of a knot measures 
its
 complexity.  It was first introduced by Bachman and Schleimer in
 \cite{BS}.  The definition we use in this paper is slightly different and corresponds more
 closely to the definition of the distance of a surface.  
 
 \begin{defin} Suppose $M$ is
 a closed, orientable irreducible 3-manifold containing a knot $K$ 
and suppose $P$ is a bridge surface for $K$ splitting $M$ into 
handlebodies $A$ and $B$. The curve complex $\mathcal{C}(P_K)$ is a graph with 
vertices corresponding to isotopy classes of essential simple closed curves on $P_K$.
Two vertices are adjacent in $\mathcal{C}(P_K)$ if their 
corresponding classes of curves have disjoint representatives.
 Let $\mcA$ (resp $\mcB$) be the set of all essential
  simple closed curves on $P_K$ that bound disks in
  $A_K$ (resp $B_K$). Then $d(P,K)=d(\mcA,\mcB)$ measured in 
  $\mathcal{C}(P_K)$. 
  \end{defin}

 The curve complex for a non-punctured torus and a 4 punctured sphere 
 are not connected. However 
 2 bridge knots in $S^3$ cannot have multiple bridge 
 surfaces, \cite{SchTo072}, so these cases don't arise in our context.   
    
 \section{Bounds on distance given by an incompressible 
surface} \label{sec:onesidebound}

 We will continue to assume that $(A,K)$ is a $K$-handlebody, $P= \bdd 
 A$ and if $A$ is a ball, then $K$ has at least 3 components.  For
 clarity we
 will refer to a properly embedded surface $E_K \subset
 A_K$ with zero Euler characteristic as an annulus only if it has 2
 boundary components both lying on $P_K$ and distinguish it from a
 punctured disk, a surface with one boundary component lying on $P_K$
 that intersects $cl(N(K))$ in a single meridional circle. Consider the 
 curve complex $\mathcal{C}(P_K)$ of $P_K$ and let $\mcA$ be the set 
of all essential 
 curves on $P_K$ that bound disks in $A_K$. 

 \begin{prop} \label{prop:cutimpliescompressing}
       Suppose $D^c$ is a cut disk for $P_K$ in $A_K$. Then there is a 
 compressing disk $D$ for $P_K$ such that $d(\bdd D^c, \bdd D)\leq 1$.   
	  
      \end{prop}

      \begin{proof}
	Let $\kappa$ be the arc of $K$ that punctures $D^c$ and $B$ be its 
	bridge disk. After perhaps an isotopy of $B$, $B\cap D^c$ 
	is a single arc $\alpha$ that separates $B$ into two 
	subdisks $B_1$ and $B_2$. Consider a regular neighborhood of 
 $D^c \cup B_1$ say. Its boundary contains a disk that intersects $P_K$ in an
	essential curve and does not intersect $\bdd D^c$ as required.  
       \end{proof}

 \begin{prop}\label{prop:essentialeasy}
    
     Consider $(F,\bdd F)\subset (A, P)$, a properly embedded surface transverse to $K$:
     \begin{itemize}
	 \item If the surface $F_K$ contains a disk component whose boundary 
	 is essential on $P_K$,
	 then $d(\mathcal{A}, f) \leq 1$ for every $f \in F_K \cap P_K$ that 
	 is essential on $P_K$.  
	 \item If $F_K$ has a punctured disk component $D^c$ whose boundary 
	 is essential on 
	 $P_K$, then $d(\bdd D^c, \mcA)\leq 1$.
	 \end{itemize}
  \end{prop}   
 \begin{proof}
  If $F_K$ contains such a disk component $D$, then $D$ is 
necessarily a 
  compressing disk for $P_K$ so
  $\bdd D \in \mathcal{A}$ and $\bdd D \cap f = \emptyset$ for every 
  $f\in F_K\cap P_K$ as $F_K$ is embedded thus
  $d(\mathcal{A}, f) \leq 1$.

  The second claim follows immediately from Proposition \ref{prop:cutimpliescompressing}. 
  \end{proof}

 \begin{prop} \label{prop:essential}
     Consider $(F,\bdd F)\subset (A, P)$, a properly embedded surface transverse to $K$
     and suppose it satisfies all of the following conditions
 \begin{enumerate}
     \item $F_K$ has no disk components,
     \item $F_K$ is c-incompressible,
     \item $F_K$ intersects every spine $\Sigma_{(A,K)}$,
\item all curves of $F_K\cap P_K$ are either essential on $P_K$ or bound punctured disks on both surfaces.

\end {enumerate}

Then there is at least one curve $f\in F_K\cap P_K$ that is essential on $P_K$ and
     such that $d(\mathcal{A}, f) \leq 1-\chi(F_K)$ and every 
       $f \in F_K\cap P_K$ that is essential on $P_K$ for which the inequality does not 
hold lies in the boundary of a $P_K$-parallel annulus component of
     $F_K$.

 \end{prop}

     \begin{proof}
If $F_K$ is a counterexample to the proposition, the surface $F_K^-$ obtained 
from $F_K$ by deleting all $P_K$-parallel annuli and $P_K$-parallel punctured disk components 
would also be a counterexample with the same euler characteristic. For $F_K^-$ is nonempty as otherwise 
$F_K$ would be disjoint from a spine $\Sigma_{(A,K)}$ and is 
c-incompressible as any c-disk would also be a c-disk for $F_K$. Thus 
we assume $F_K$ does no have any $P_K$-parallel annuli or punctured disk components.

Let $E$ be a compressing disk for $P_K$ in $A_K$ (not punctured by the
knot) so that $|E \cap F_K|$ is minimal among all such disks.  If in
fact $E\cap F_K = \emptyset$, then $d(\bdd E, f) \leq 1$ for every $f
\in \bdd F_K$ as required so we may assume $E\cap F_K \neq \emptyset$.  Circles of intersection between $F_K$
and $E$ and arcs that are inessential on $F_K$ can be removed by
innermost disk and outermost arc arguments.  Thus we can assume $F_K$
and $E$ only intersect in arcs that are essential on $F_K$.

The proof now is by induction on $1-\chi(F_K)$.  As $F_K$ has no disk
components for the base case of the induction assume $1-\chi(F_K) =
1$, i.e. all components of $F_K$ are annuli or once punctured disks 
and
no component is $P_K$-parallel.  If $E$
intersects a punctured disk component of $F_K$ the arc of intersection
would necessarily be inessential on $F_K$ contradicting the minimality
of $|F_K \cap E|$ so we may assume that if $F_K \cap E \neq 
\emptyset$, $E$ only intersects
annulus components of $F_K$.  An outermost arc of intersection on $E$ 
bounds a 
$P$-compressing disk $E_0$ for $F_K$.  After the $P$-compression, the
new surface $F_K'$ contains a compressing disk $D$ for $P_K$, the
result of a $P$-compression of an essential annulus, and $\bdd D$ is
disjoint from all $f \in \bdd F_K$ by Lemma 
\ref{lemma:boundredessential}.  As $\bdd D \in \mathcal{A}$,
$d(f, \mathcal{A})\leq 1= 1-\chi(F_K)$ for every $f\in F_K\cap P_K$ 
as desired.

Now suppose $1-\chi(F_K) > 1$. 
Again let $E_0$ be a subdisk of $E$ 
cut off by an outermost arc of $E\cap F_K$ and $F_K'$ be the surface 
obtained after the $P$-compression.  By Lemma 
\ref{lemma:boundredessential} $F_K'$ also intersect every spine 
$\Sigma_{(A,K)}$ and is c-incompressible.
By the definition of $P$-compression, $F_K'$ cannot have any disk 
components as $F_K$ did not have any. Thus $F_K'$ satisfies the first 
3 conditions of the proposition. There are two cases to consider.

{\bf Case 1:} Any simple closed curves in $F_K'\cap P_K$ that are inessential on $P_K$
bound punctured disks on both surface. 

In this case $F_K'$ satisfies all the hypothesis of the proposition so we 
can apply the induction hypothesis. Thus
there exists a curve $f'\in F_K' \cap P_K$ that satisfies the distance inequality. 
Since, by Lemma \ref{lemma:boundredessential}, for every component $f$ 
of 
$F_K\cap P_K$, $d(f, f') \leq 1$, we have
the inequality $d(f, \mathcal{A}) \leq d(f',\mathcal{A}) + d(f, f')
\leq 1-\chi(F_K')+1=1-\chi(F_K)$, as desired.

{\bf Case 2:} Some curve of $F_K' \cap P_K$ is inessential on $P_K$ 
but does not bound a punctured disk on $F_K'$.

Let $c$ be this curve and let $E^*$ be the possibly punctured disk 
$c$ bounds on $P_K$. By 
our hypothesis, the tunnel dual to the 
$P$-compression must be adjacent to $c$ as otherwise $c$ would 
persist in $F_K \cap P_K$. Push a copy of $E^*$ slightly into $A_K$. 
After the tunneling, $E^*$ is no longer parallel to $P_K$. As 
$F_K$ was assumed to be c-incompressible, $c=\bdd E^*$ must be parallel 
to 
some component of $\bdd F_K$. As $c$ didn't bound a punctured 
disk on $F_K'$, $\bdd E^*$ must be parallel to some component $\tilde 
c\in F_K\cap P_K$ that is essential on $P_K$ by hypothesis. 
Use this parallelism to extend $E^*$ to a c-disk 
for $P_K$ with boundary $\tilde c$, see Figure \ref{fig:extension}. Now 
for every $f \in F_K \cap P_K$, by Proposition
\ref{prop:essentialeasy} we have that $d(f, 
\mcA)\leq d(f, \bdd E^*)+ d(\bdd E^*,\mcA)\leq 1+1=2\leq1-\chi(F_K)$. 

\begin{figure}[tbh]
\centering
 \includegraphics[width=.5\textwidth]{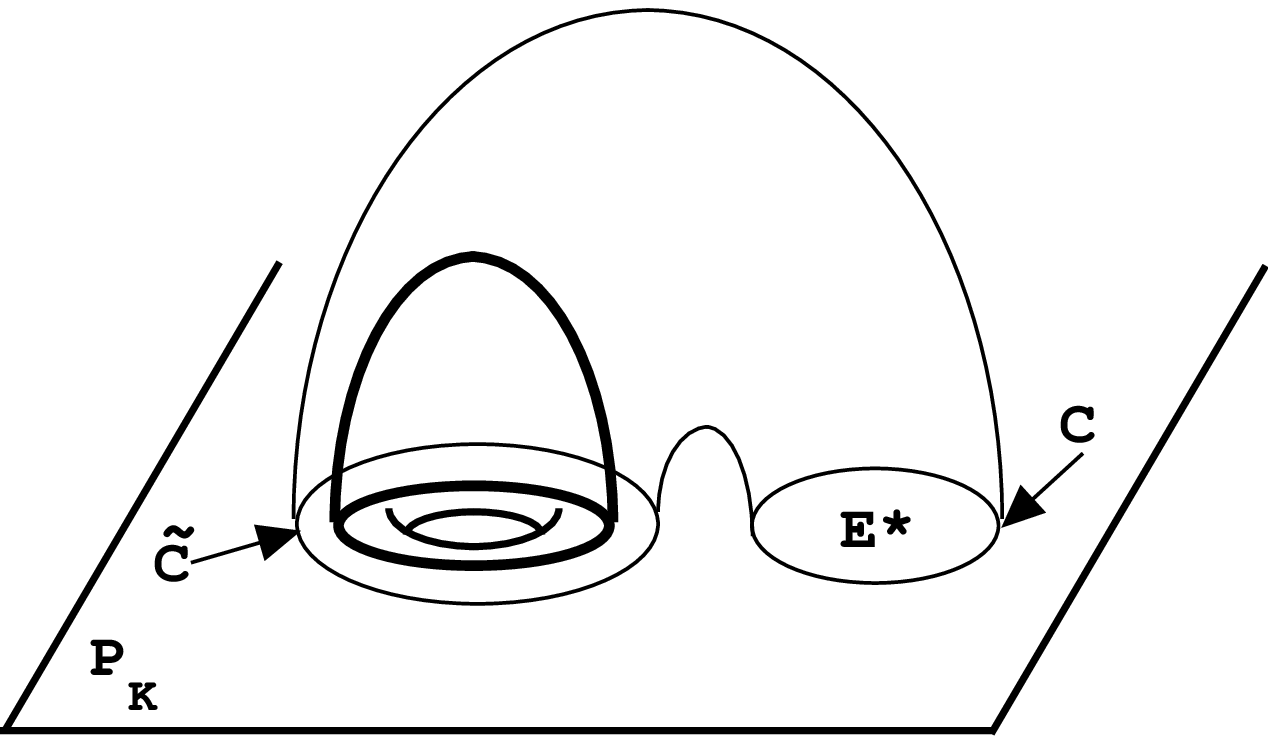}
    \caption{} \label{fig:extension}
   \end{figure}
\end{proof}

\section{The genus of an essential surface bounds the distance
       of a knot}\label{sec:essentialbound}

\begin{notation} For the rest of the paper we will assume that $M$ is
a closed irreducible orientable 3-manifold containing a knot $K$ and $P$ is 
a bridge surface for $K$ such that $M=A\cup_P B$.  Furthermore we 
assume that if $P$ is a sphere, then $P_K$ has at least 6 punctures.
\end{notation}

Let $Q \subset M$ be a properly embedded surface that is transverse to
$K$.  We will consider how the surfaces $P_K$ and $Q_K$ can intersect
in $M_K$ to obtain bounds on $d(P,K)$. 

We import the next lemma directly from \cite{ST}.

     \begin{lemma} \label{lem:essentialdist}
	 Let $Q \subset M$ be a properly embedded surface that is transverse to
	 $K$ and let $Q_{K}^{A} = Q_K \cap A_K,
	 Q_K^{B} = Q_K \cap B_K$. Suppose $Q_K$ satisfies the following
    conditions:
    \begin{itemize}
     \item All curves of $P_K \cap Q_K$ are essential on $P_K$ and
     don't bound disks on $Q_K$.

     \item
     There is at least one curve $a \in Q_K^A \cap P_K$ such that
     $d(a,\mathcal{A}) \leq 1-\chi(Q_K^A)$ and any curve in $Q_K^A \cap
     P_K$ for which the inequality does not hold is the boundary of an
     annulus component of $Q_K^A$ that is parallel into $P_K$.

\item There is at least one curve $b \in Q_K^B \cap P_K$ such that
$d(b,\mathcal{B}) \leq 1-\chi(Q_K^B)$ and any curve in $Q_K^B
\cap P_K$ for which the inequality does not hold is the
boundary of an annulus component of $Q_K^B$ that is parallel
into $P_K$.
    \end{itemize}

Then $d(K,P) \leq 2-\chi(Q_K)$ .

     \end{lemma}
     \begin{proof}
     Call a component $c$ of $P_K \cap Q_K$ {\em A-conforming} (resp
     {\em B-conforming}) if $d(c,\mathcal{A}) \leq 1-\chi(Q_K^A)$
     (resp $d(c,\mathcal{B}) \leq 1-\chi(Q_K^B)$).  By hypothesis
     there are both A-conforming components of $Q_K \cap P_K$ and
     B-conforming components.  If there is a component $c$ that is
     both A-conforming and B-conforming, then $$d(K,P) =
     d(\mathcal{A}, \mathcal{B}) \leq d(\mathcal{A},c) +
     d(c,\mathcal{B}) \leq 2 - \chi(Q_K^{A}) - \chi(Q_K^{B}) = 2 -
     \chi(Q_K)$$ as required.

     If there is no such component, let $\gamma$ be a path in $Q_K$
     from an A-conforming component to a B-conforming component,
     chosen to intersect $P_K$ as few times as possible.  In
     particular, any component of $P_K \cap Q_K$ incident to the
     interior of $\gamma$ is neither A-conforming nor B-conforming, so
     each of these components of $Q_K^{A}$ and $Q_K^{B}$ is an
     annulus, parallel to an annulus in $P_K$.  It follows that the
     components of $P_K\cap Q_K$ at the ends of $\gamma$ are isotopic
     in $P_K$.  Letting $c$ be a simple closed curve in that isotopy
     class in $P_K$ we have as above $$d(K,P) = d(\mathcal{A},
     \mathcal{B}) \leq d(\mathcal{A},c) + d(c,\mathcal{B}) \leq 2 -
     \chi(Q_K^{A}) - \chi(Q_K^{B}) = 2 - \chi(Q_K)$$ as required. 
     \end{proof}

     \begin{cor} \label{cor:essentialdist}

Suppose $Q_K \subset M_K$ is a properly embedded connected surface
transverse to $P_K$ so that all curves of $P_K \cap Q_K$ are essential
on both surfaces.  If $Q_K^{A}$ and $Q_K^{B}$ are c-incompressible and
intersect every spine $\Sigma_{(A,K)}$ and $\Sigma_{(B,K)}$
respectively , then $d(K,P) \leq 2-\chi(Q_K)$ .

     \end{cor}

	 \begin{proof} Proposition \ref{prop:essential} shows that
	 $Q_K^{A}$ and $Q_K^{B}$ satisfy respectively the second and
	 third conditions of Lemma \ref{lem:essentialdist}.
     \end{proof}

      The following definition was first used by Scharlemann in
      \cite{ST}.
      
      \begin{defin} \label{def:removable} Suppose $S$ and $T$ are two 
      properly embedded surfaces in a 3-manifold $M$ containing a 
      knot $K$ and assume $S$ and $T$ intersect the knot 
      transversely. Let $c\in S_K\cap T_K$
     be a simple closed curve bounding possibly punctured disks
      $D\subset S_K$ and $E \subset T_K$.  If $D$ intersects $T_K$ only in curves that are inessential in
      $T_K$ and
      $E$ intersects $S_K$ only in curves that are inessential in
      $S_K$ we say that $c$ is {\em removable}.
	       
      \end{defin}
      The term reflects the fact that all such curves can be removed
      by isotopies of $S_K$ whose support lies away from any curves of
      intersection that are essential either in $S_K$ or in $T_K$. 
      Indeed, if $c$ is removable, then any component of $D \cap E$ is
      clearly also removable.
      
      The following definition was introduced by Bachman and Schleimer in \cite{BS}.

     \begin{defin}
      Suppose $S$ and $T$ are two properly embedded surfaces in a
     3-manifold $M$.  A simple closed curve
     $\alpha \in S\cap T$ is mutually essential if it
     is essential on both surfaces, 
     it is mutually inessential if it is inessential on both surfaces and it is mutual if
     it is either mutually essential or mutually inessential.
		  
	       \end{defin}
      The following remark follows directly from the above two definitions. 
      \begin{rmk}\label{rmk:mutualremovable}
	If every curve of intersection between $S_K$ and $T_K$ is mutual, then all 
inessential curves of $S_K \cap T_K$ are removable.   
\end{rmk}	  
      
Now we can recover the bound on distance obtained in \cite{BS} but
using our definition of distance.  Note that we only require the 
surface $Q_K$ to
have no compressing disks but allow it to have cut-disks.

     \begin{thm}\label{thm:essbound}
    Let $M$ be a closed irreducible orientable manifold 
     containing a knot $K$ and let $P$ be a bridge surface for $K$ 
     such that if $P$ is a sphere, $P_K$ has at least 6 punctures. 
     Suppose $Q \subset M$ is a properly embedded essential (in $M_K$)
     meridional surface such that $Q_K$ is neither a sphere nor an annulus.  Then $d(K,P) \leq
     2-\chi(Q_K)$.  If $Q_K$ is an essential annulus, then $d(K,P) \leq 3$
     \end{thm}

     \begin{proof}
	 
If $Q_K$ has any cut disks, cut-compress along them, i.e. if $D^c$ is
a cut disk for $Q_K$, remove a neighborhood of $\bdd D^c$ from $Q_K$
and then add two copies of $D^c$ along the two newly created boundary
components.  Repeat
this process until the resulting surface has no c-disks. Let $Q_K'$ be the resulting surface and notice that 
$\chi(Q_K)=\chi(Q_K')$. 
Suppose $Q_K'$ has a compressing disk $D$. 
$Q_K$ can be recovered from $Q_K'$ by tubing along a collection of 
subarcs of $K$. Note that as $D \cap K=\emptyset$ none of these tubes 
can intersect $D$. Thus $D$ is also a compressing disk for $Q_K$ 
contrary to the hypothesis so $Q_K'$ is also incompressible. Finally note 
that in this
process no sphere, annulus or torus components are produced so at least one of the
resulting components is not a sphere, annulus or torus, in particular $Q_K'$ has 
at least one component that is not parallel to $\bdd M_K$. By 
possibly replacing $Q_K$ by $Q_K'$ we may assume that $Q_K$ is also 
cut-incompressible.

 Recall that $\Sigma_{(A,K)}$ and $\Sigma_{(B,K)}$ are the spines for
 the $K$-handlebodies $(A,K)$ and $(B,K)$.  Consider $H: P_K \times
 (I, \bdd I) \rightarrow (M_K, \Sigma_{(A,K)} \cup \Sigma_{(B,K)})$, a
 sweep-out of $P_K$ between the two spines.  For a fixed generic value
 of $t$, $H(P_K, t)$ will be denoted by $P_K^t$.  By slightly abusing
 notation we will continue to denote by $A_K$ and $B_K$ the two
 components of $M_K-P_K^t$ and let $Q_K^A=Q_K\cap A_K$ and
 $Q_K^B=Q_K\cap B_K$.
  
During the sweep-out, $P_K^t$ and $Q_K$ intersect generically except
in a finite collection of values of $t$.  Let $t_1,..t_{n-1}$ be these
critical values separating the unit interval into regions where
$P_K^t$ and $Q_K$ intersect transversely.  For a generic value $t$ of
$H$, the surfaces $Q_K$ and $P_K^t$ intersect in a collection of
simple closed curves.  After removing all removable curves, label a
region $(t_i, t_{i+1}) \subset I$ with the letter $A^*$ (resp $B^*$) if
$Q_K^A$ (resp $Q_K^B$) has a disk or punctured disk component in the region whose
boundary is essential on $P_K$.  

Suppose $Q_K^A$ say, can be isotoped off some spine $\Sigma_{(A,K)}$. 
Then, using the product structure between the spines and the fact that
all boundary components of $Q_K$ lying on the knot are meridional, we
can push $Q_K$ to lie entirely in $B_K$ contradicting Lemma \ref
{lem:notclosed}.  Therefore $Q_K$ must intersect both spines
$\Sigma_{(A,K)}$ and $\Sigma_{(B,K)}$ in meridional circles and so the
subintervals adjacent to the two endpoints of the interval are labeled
$A^*$ and $B^*$ respectively.

{\bf Case 1:} Suppose there is an unlabeled region.  If some curve of
$Q_K \cap P_K$ is inessential on $P_K$ in that region, it must also be
inessential on $Q_K$ as otherwise it would bound a c-disk for $Q_K$. 
Suppose some curve is essential on $P_K$ but inessential on $Q_K$. 
This curve would give rise to one of the labels $A^*$ or $B^*$
contradicting our assumption.  We conclude that all curves of $P_K
\cap Q_K$ are mutual.  In fact this implies that all curves $P_K\cap Q_K$ are
essential on $Q_K$ and on $P_K$ as otherwise they would be removable 
by Lemma \ref{lem:removable} and all removable curves have already 
been removed. 
Suppose $Q_K^A$ say has a c-disk.  The boundary of this c-disk would
also be essential in $Q_K$ contradicting the hypothesis thus we
conclude that in this region $Q_K^A$ and $Q_K^B$ satisfy the
hypothesis of Corollary \ref{cor:essentialdist} and thus $d(K,P)\leq
2-\chi (Q_K)$.

{\bf Case 2:} Suppose there are two adjacent regions labeled $A^*$ 
and $B^*$.  (This includes the case when one or both of
these regions actually have both labels)

The labels are coming from possibly punctured disk components of
$Q_K-P_K$ that we will denote by $D^*_A$ and $D^*_B$ respectively. 
Using the triangle inequality we obtain
\begin{equation}\label{eq1}
    d(K,P) \leq d(\mcA, \bdd D^*_A) + d(\bdd D^*_A, \bdd D^*_B)+d(\bdd
    D^*_B, \mcB).
\end{equation}

The curves of intersection before and after going through the critical
point separating the two regions can be made disjoint so $d(\bdd
D^*_A, \bdd D^*_B)\leq 1$ (the proof of this fact is similar to the
proof of the last item of Lemma \ref{lemma:boundredessential}).  By
Proposition \ref{prop:essentialeasy} $d(\mcA, \bdd D^*_A),d(\mcB, \bdd
D^*_B)\leq 1$ so the equation above gives us that $d(K,P)\leq 3\leq
2-\chi (Q_K)$ as long as $\chi(Q_K)<0$.

If $\chi(Q_K)= 0$ and $Q_K$ is a torus, $D^*_A$ and $D^*_B$ must be
disks, so $d(\mcA, \bdd D^*_A)=d(\mcB, \bdd D^*_B)=0$.  Thus Equation
\ref{eq1} gives us that $d(K,P)\leq d(\bdd D^*_A, \bdd D^*_B)\leq
1\leq 2-\chi (Q_K)$.  If $Q_K$ is an essential annulus, we conclude that
$d(K,P)\leq 3$

     \end{proof}

     \begin{cor}
Suppose $K=K_1 \# K_2$, then any bridge surface for $K$ has distance 
at most 3.
     
     \end{cor}
     
     \begin{proof}
The sphere that decomposes $K$ into its factors is an essential 
annulus in $M_K$.	 
\end{proof}	 
	 
     \section{Edgeslides} \label{sec:edgeslides}

This section is meant to provide a brief overview of edgeslides as 
first described in \cite{RS}. Here we only give sketches of 
the relevant proofs and references for the complete proofs.   

Suppose $(Q, \bdd Q) \subset (M, P)$ is a bicompressible splitting
surface in an irreducible 3-manifold with $P \subset \bdd M$ a 
compact 
sub-surface, (in our context $M$ will be a $K$-handlebody
 and $P$ its punctured 
boundary).  Let $X, Y$ be the two
components of $M-Q$ and let $Q_X$ be the result of maximally
compressing $Q$ into $X$. The compressions can be undone by 
tubing along the edges of a graph
$\Gamma$ dual to the compressing disks.  We will denote by $X^-$ and 
$Y^+$ the components of $M-Q_X$ with $X
\supset X^-$ and $Y \subset Y^+$, in particular $\Gamma \subset Y^+$. 
Let $\Delta \subset Y$ be a set of
compressing disks for $Q$, thus $\bdd \Delta \subset Q_X \cup
\Gamma$.  Finally $T$ will be a disk in $Y^+$ with $\bdd T \subset 
(Q_X \cup
P)$ that is not parallel to a subdisk of $Q_X \cup P$ and $\Lambda$ 
will be the graph on $T$ with vertices $T \cap \Gamma$ and
edges $T \cap \Delta$. 

The graph $\Gamma$ described above is not unique, choosing a 
different 
graph is equivalent to an isotopy of $Q$. All graphs that are dual to 
the same set of compressing disks are related by {\em edge slides}, 
i.e. 
sliding the endpoint of some edge along other edges of $\Gamma$. The 
precise definition can be found in \cite{Saito} or \cite {STh}.
	
The following lemma is quite technical, a detailed proof of a 
very similar result can be found in \cite{Saito} Proposition
3.2.2 or \cite {STh}, Prop.  2.2.  We will only briefly sketch the
proof here but we will provide detailed references to the 
corresponding results in \cite{Saito} and note that there the letter 
$P$ is used for
the disk we call $T$ but all other notation is identical.

	\begin{lemma}
 \label{lem:edgeslides}
	  If $\Lambda$ cannot have isolated vertices,
	then we can take $\Gamma$ and $\Delta$ to be disjoint from $T$ by
	\begin{itemize}
	\item isotopies of $T$ rel $\bdd T$
	\item rechoosing $\Delta$ keeping
	$|\Delta|$ fixed. 
	\item edge slides of $\Gamma$ (which translate into isotopies of 
$Q$).
	\end{itemize}

	       \end{lemma}

	    \begin{proof}

	 Pick an isotopy class of $T$ rel. $\bdd T$, an isotopy class of
	 $\Delta$ and a representation of $\Gamma$ such that $(|T \cap
	 \Gamma|, |T \cap \Delta|)$ is minimal in the lexicographic order.

	{\bf Claim 1:} Each component of $T \cap \Delta$ is an arc
	(\cite{Saito}, Lemma 3.2.3).

	Suppose $T \cap \Delta$ contains a closed curve component.  The
	innermost such on $\Delta$, $\omega$ bounds a disk $D_0$ on $\Delta$
	disjoint from $T$.  Via an isotopy of the interior of $T$, using the
	fact that $M$ is irreducible, the disk
	$\omega$ bounds on $T$ can be replaced with $D_0$ thus eliminating at
	least $\omega$ from $T \cap \Delta$ contradicting minimality.

	{\bf Claim 2:} $\Lambda$ has no inessential loops, that is loops 
that 
	bound disks in $T-\Gamma$ (\cite{Saito} Lemma
	3.2.4).

	Suppose $\mu$ is a loop in $\Lambda$ and let $D \in \Delta$ be such
	that $\mu \subset D$.  The loop $\mu$ cuts off a disk $E \subset T$.
	As a subset of $D$, $\mu$ is an arc dividing $D$ into two subdisks
	$D_1$ and $D_2$.  (The disk $E$ resembles a boundary compressing disk
	for $D$ if we think of $\eta (\Gamma)$ as a boundary component.)
	At
	least one of $D_1 \cup E$ and $D_2 \cup E$ must be a compressing 
disk.  Replace $D$ with this disk
	 reducing $|T \cap \Delta|$. 
	 
	 {\bf Claim 3:} $\Lambda$ has no isolated vertices (\cite{Saito}, 
Lemma 3.2.5).

		 By hypothesis.
		 
	{\bf Claim 4:} Every vertex of $\Lambda$ is a base of a loop
	(\cite{Saito}, Lemma 3.2.6).

	Suppose $w$ a vertex of $\Lambda$ is not a base of any loop, we will
	show we can reduce $(|T \cap \Gamma|, |T \cap \Delta|)$.

	Let $\sigma$ be the edge of $\Gamma$ such that $w \in \sigma \cap T$.
	As $w$ is not isolated, there is a disk $D \in \Delta$ such that $w
	\in \bdd D$.
	$D \cap T$ is a collection of arcs that are edges in $\Lambda$. Let
	$\gamma$ be an
	outermost arc on $D$ of all arcs that have $w$ as one endpoint. Let
	$w'$ be the other end point of $\gamma$. Then $\gamma$ cuts a subdisk
	$D_{\gamma}$ from $D$ the interior of which may intersect $T$ but
	$\bdd D_{\gamma}$ only contains one copy of $w\in \bdd \gamma$. Thus
	there cannot be an entire copy of the edge $\sigma$ in $\bdd
	D_{\gamma}$ and so there are three possibilities.

	Case 1: $(\bdd D_{\gamma} - \gamma) \subset \sigma$.  Then we can
	perform an edge slide of $\sigma$ which removes $\gamma$ from
	$\Lambda$. (\cite{Saito} Fig. 23).

	Case 2: $(\bdd D_{\gamma} - \gamma)$ contains some subset of
	$\sigma$ with only one copy of one of the
	endpoints
	of $\sigma$. By sliding  $\sigma$ along
	$D_{\gamma}$ we can reduce this case to the first case. 
(\cite{Saito} Fig.
	24)

	Case 3: $(\bdd D_{\gamma} - \gamma)$ contains some subset of
	$\sigma$ but it contains two copies of the same endpoint of 
$\sigma$. This
	is the most complicated case requiring
	{\em broken edge slides} and (\cite {Saito}, Fig. 25)  has an 
excellent discussion
	on the topic. 

	By the above 4 claims we can conclude that $\Lambda= \emptyset$ as 
	desired, for by claim 4 some loop must be inessential contradicting 
	claim 2.

	\end{proof}

	\begin{rmk}
	If $Q$ is weakly incompressible, the hypothesis of the lemma are 
satisfied 
	as a meridional circle of an isolated vertex of $\Lambda$ will be a 
	compressing disk for $Q$ in $X$ that is disjoint from the set of 
	compressing disks $\Delta \in Y$.
	\end{rmk}

	    \begin{cor}  \label{cor:maxincomp}

	Let $(Q, \bdd Q) \subset (M, \bdd M)$ be a bicompressible weakly 
incompressible
	surface splitting $M$ into component $X$ and $Y$.  Let $Q_{X}$ be 
the result
	of maximally compressing $Q$ into $X$.  Then $Q_{X}$ is 
incompressible in $M$.
	    \end{cor}

	\begin{proof}
	  The argument is virtually identical to the argument in \cite{Sc2}.
	  Suppose $Q_X$ is
	  compressible with compressing disk $D$ that necessarily lies in
	  $Y^+$. Let $E$ be a compressing disk for $Q$ in $Y$. As $Q$ is
	  weakly incompressible, by the above remark we can apply Lemma
	\ref{lem:edgeslides}, with $D$
	  playing the role of $T$, and $\Delta=E$. By Lemma \ref{lem:edgeslides} we can 
arrange
	that $(E \cup \Gamma) \cap D
	=\emptyset$
	  so $D$ is also a compressing disk for $Q$ in $Y$ and is disjoint
	  from $\Gamma$ and thus from all compressing disks for $Q$ in $X$ 
	  contradicting weak incompressibility of $Q$.
       \end{proof}

\begin{cor}\label{cor:nicecompdisk}
       Suppose $(A,K)$ is a $K$-handlebody with $\bdd A=P$ and $F$ is 
a 
       bicompressible surface splitting $A$ into submanifolds $X$ and 
       $Y$. Let $F_K^X$ be the result of maximally compressing $F_K$ 
       into $X_K$. Then there exists a compressing disk $D$ for $P_K$ 
       that is disjoint from a complete collection of compressing 
       disks for $F_K$ in $X_K$ and intersects $F_K$ only in arcs 
       that are essential on $F_K^X$. 
       
       \end{cor}
       
       \begin{proof}
	   
	   Select a disk $D$ and isotope $F_K^X$ to minimize $|D\cap F_K^X|$ 
	   and 
	   then choose a representation of $\Gamma$ that minimizes $|D\cap 
	   \Gamma|$. As $A-N(K)$ is irreducible, by an innermost disk and 
	   outermost arc arguments, $D$ intersect $F_K^X$ in essential arcs 
	   only. Applying Lemma \ref{lem:edgeslides} with the 
	   disk $T$ playing the role of $D$, we conclude that $\Gamma$ is 
	   disjoint from $D$. 
\end{proof}

\section{Bounds on distance given by a c-weakly incompressible 
surface} \label{sec:keytheorem}

  Our ultimate goal in this paper is to extend Theorem 
\ref{thm:essbound} to allow for
both $P$ and $Q$ to be bridge surfaces for the same knot.  To do this, we need a theorem similar to
 Proposition \ref{prop:essential} but allowing for $F_K$ to have 
 certain kinds of c-disks.

 \begin{notation}
     In this section let $(A,K)$ be a $K$-handlebody with boundary 
     $P$ such that if $A$ is a ball, $K$ has at least 3 components and let
  $F \subset A$ be a properly embedded surface transverse to
  $K$ splitting $A$ into submanifolds $X$ and $Y$. 
\end{notation}

  \begin{defin}
     The surface $F_K$ associated to $F$ is called bicompressible 
    if $F_K$ has some compressing disks in both $X_K$ and $Y_K$. The 
surface is called 
    cut-bicompressible if it has cut-disks in both $X_K$ and $Y_K$. 
    Finally, the surface is called c-bicompressible if it has c-disks 
    in both $X_K$ and $Y_K$.
  \end{defin}

   The next definition is an adaptation of the idea of a weakly
  incompressible surface but taking into consideration not only
  compressing disks but also cut disks.
  
\begin{defin}
  The surface $F_K$ is called c-weakly incompressible if it is
  c-bicompressible and any pair $D^*_X, D^*_Y$ of c-disks contained 
  in $X_K$ and $Y_K$ respectively intersect along their boundary.
 \end{defin}

\begin{prop} \label{prop:nosinglepoint}
  If a splitting surface $F_K \subset A_K$ has a 
    pair of two compressing disks or a compressing disk and a cut disk that are 
    on opposite sides of $F_K$ and intersect in exactly one point, 
then $F_K$ is c-strongly compressible.    

\end{prop}
 \begin{proof}
     Suppose $F \subset A$ splits $A$ into manifolds $X$ and $Y$ and 
let $D_X \subset X$ and $D_Y\subset Y$ be a pair of disks that intersect in exactly one point. 
 Then a neighborhood of $D_X\cup D_Y$ contains a pair of compressing disks on 
 opposite sides of $F_K$ with disjoint boundaries (in fact their 
boundaries are isotopic). If $D_X$ say is a compressing disk 
 and $D_Y$ is a cut disk, banding two copies $D_X$ together along $\bdd D_Y$ 
 produces a compressing disk disjoint from $D_Y$, see Figure \ref{fig:recovery}.

\begin{figure}[tbh]
\centering
 \includegraphics[width=.5\textwidth]{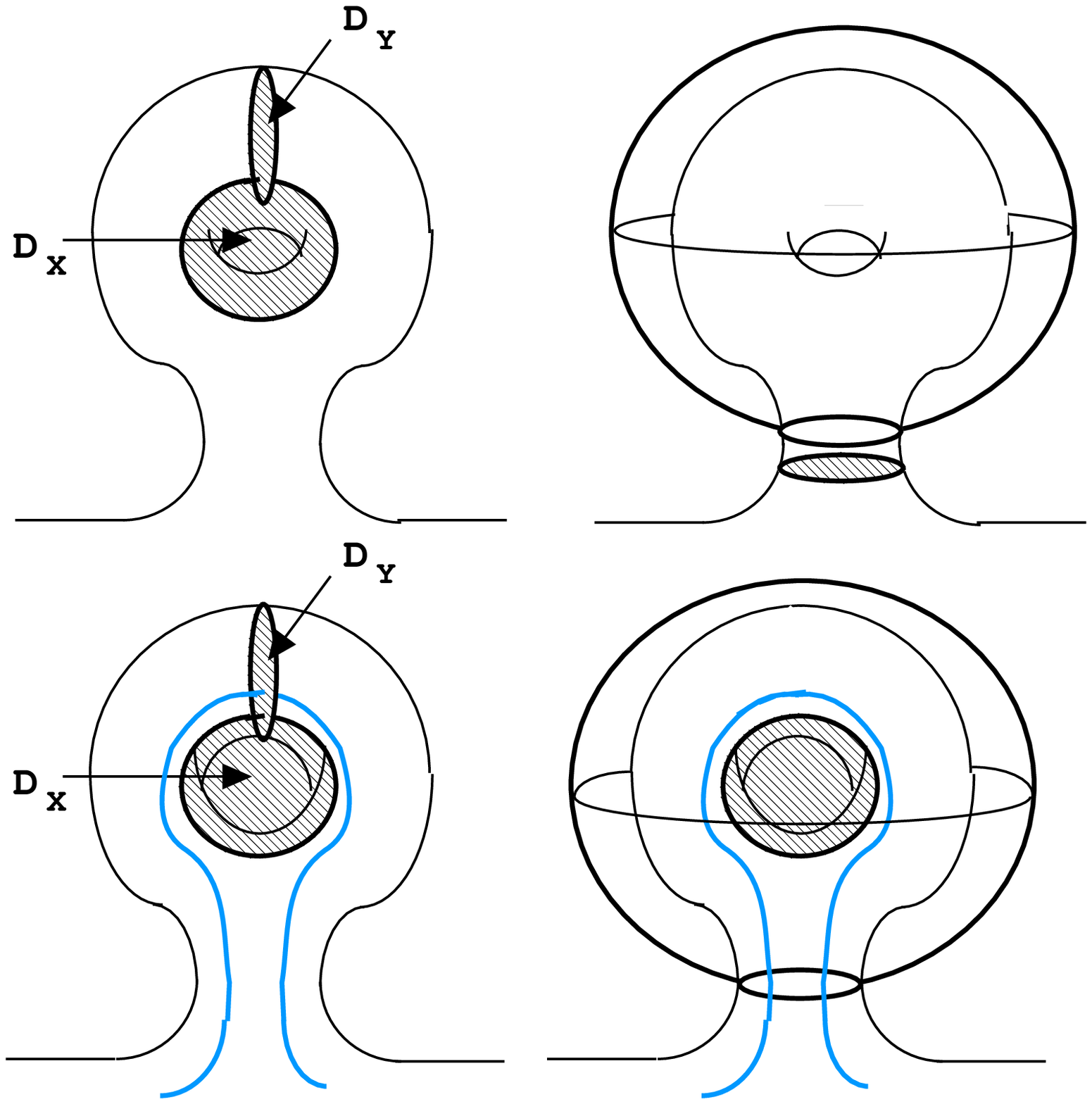}
    \caption{} \label{fig:recovery}
   \end{figure}
   
\end{proof}

   \begin{prop} \label{prop:preserving}
    Let $F_K\subset A_K$ be a c-weakly incompressible splitting surface such 
that every component of $F_K \cap P_K$ is mutual and let $F_K'$ be the
   surface obtained from $F_K$ via a $P$-compression.  If $F_K'$ is
   also c-bicompressible, then every component of $F_K' \cap
   P_K$ is essential on $P_K$ or is mutually inessential.
   \end{prop}
 
   \begin{proof}
    Let $X$ and $Y$ be the two components of $A-F$. 
    Without loss of generality, let $E_0 \subset X_K$ be the $P$-compressing 
disk for $F_K$. Suppose that there is some $f' \in
       \bdd F_K'$ that bounds a possibly punctured disk $D_{f'}$ on
       $P_K$ but not on $F_K'$.  $F_K$ can be recovered by tunneling 
       $F_K'$ along an arc $e_0 \subset P_K$. As all curves of $F_K 
       \cap P_K$ are mutual, $e_0 \cap f' \neq \emptyset$. 
       
      \textbf{Case 1:}  $e_0$ has one boundary component on $f'$ and the other on 
       some other curve $c \in P_K \cap F_K$ ($c$ may or may not be 
       essential on $P_K$). If $c \subset D_f'$, then $F_K \cap P_K$ also has a 
       curve that is inessential on $P_K$ but essential on $F_K$ 
       contrary to the hypothesis. If $c$ does not lie in $D_f'$ 
       then by slightly pushing the disk $D=D_f' \cup E_0$ away from 
       $P_K$ we obtain a c-disk for $F_K$ contained in $X_K$, see Figure 
       \ref{fig:preservinga}. By 
       hypothesis $F_K'$ is c-bicompressible, in 
       particular there is a c-disk $D'$ for $F_K'$ that lies on the 
       other side of $F_K'$ then the side $D_f'$ lies on. $D'$ is 
       also a c-disk for $F_K$ lying in $Y$ that is disjoint from 
       $D\subset X$ contradicting 
       the c-weak incompressability of $F_K$.
       \begin{figure}[tbh]
	    \centering
	    \includegraphics[width=.7\textwidth]{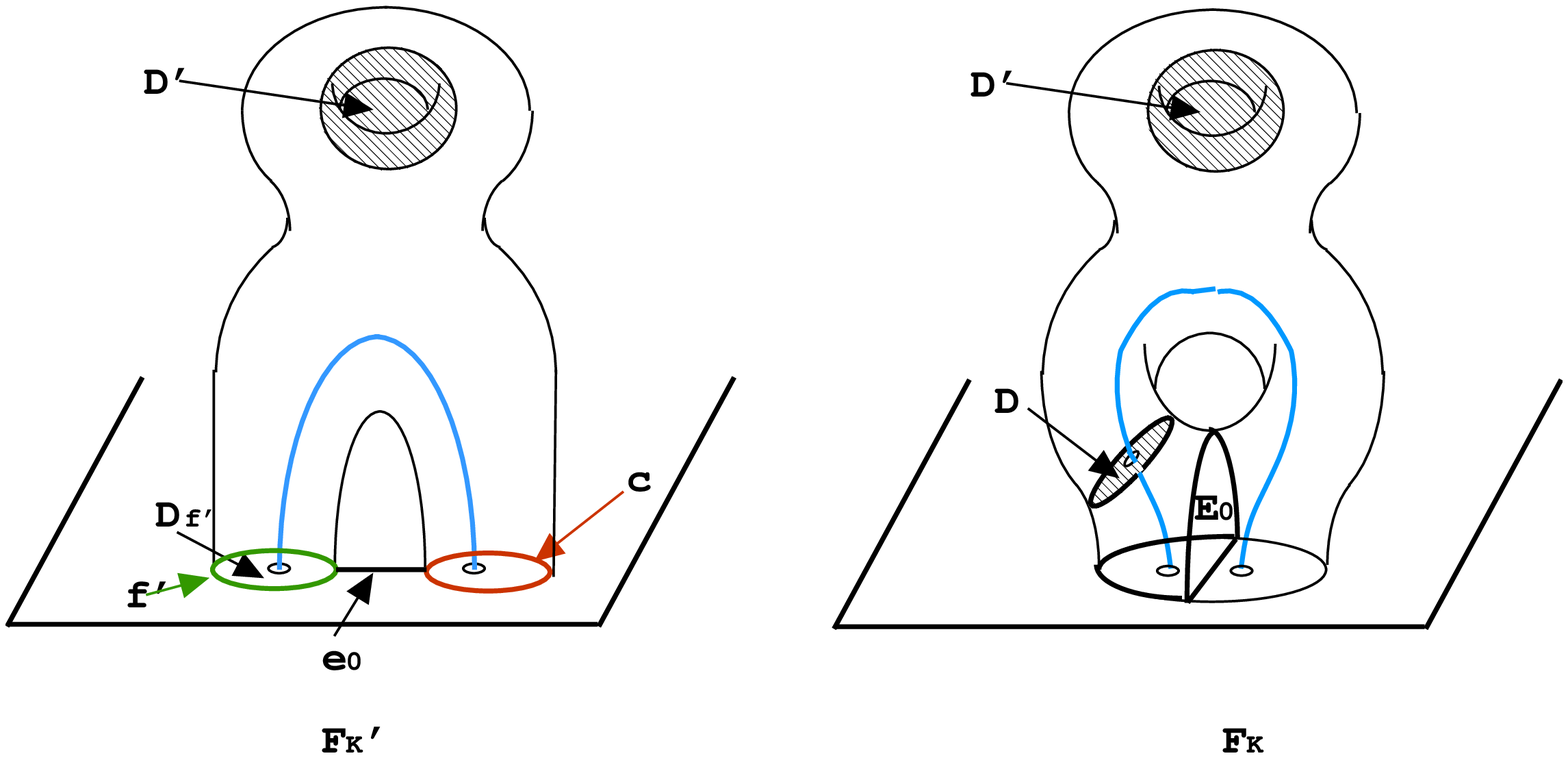}
	    \caption{} \label{fig:preservinga}
	     \end{figure}

       \textbf{Case 2:}  $e_0$ has both boundary component on $f'$. 
       If $e_0\subset D_f'$ then again $F_K \cap P_K$ has a 
       curve that is inessential on $P_K$ but essential on $F_K$ 
       contrary to the hypothesis so assume $e_0 \cap D_f'=\bdd e_0$, 
       see Figure \ref{fig:preservingb}.
       Consider the possibly punctured disk $D$ obtained by taking the union of $D_{f'}$
       together with two copies of $E_0$.
      As in the previous case this is a c-disk for $F_K$ lying in $X_K$ that is disjoint from
at least one c-disk for $F_K$ lying in $Y_K$ contradicting c-weak 
incompressibility of $F_K$.

    \begin{figure}[tbh]
      \centering
      \includegraphics[width=.8\textwidth]{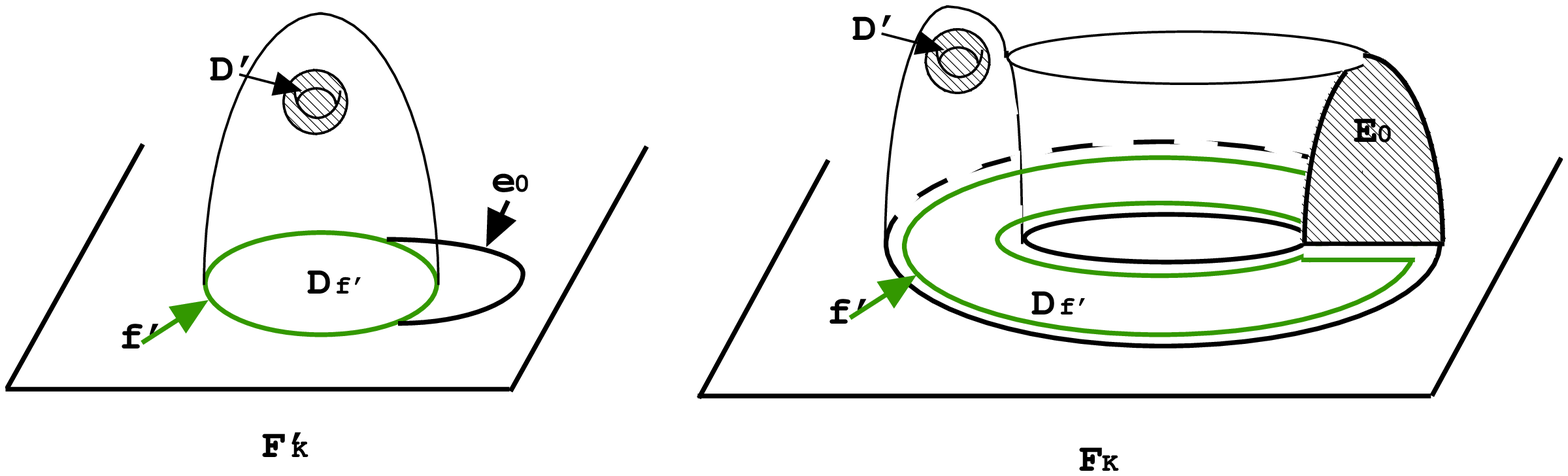}
      \caption{} \label{fig:preservingb}
       \end{figure}

   \end{proof}

     \begin{prop}\label{prop:recovery}
Suppose $F_K'$ splitting $A_K$ into $X_K'$ and $Y_K'$ satisfies one of
the following two conditions
 
	 \begin{itemize}
	     \item There is a spine $\Sigma_{(A,K)}$ entirely
	     contained in $X_K'$ say and $F_K'$ has a c-disk in $X_K'$
	     disjoint from that spine, or
	     
	     \item There is at least one curve $f' \in F_K'\cap P_K$
	     that is essential on $P_K$ and $d(f', \mcA) \leq
	     1-\chi(F_K')$.
	    
	    \end{itemize}
	    
 If $F_K$ is obtained from $F_K'$ by tunneling or tubing (possibly
 along subarcs of $K$) with all tubes lying in $Y_K'$, then $F_K$
 satisfies one of the following conditions.
 \begin{itemize}
 \item There is a spine $\Sigma_{(A,K)}$ entirely contained in $X_K$,
 and $F_K$ has a c-disk in $X_K$ disjoint from that spine, or
	      
	      \item For every curve $f$ in $F_K\cap P_K$ that is 
	      essential on $P_K$ the inequality $d(f,
	      \mcA) \leq 1-\chi(F_K)$ holds.

	      \end{itemize}
 \end{prop}
 
 \begin{proof}
 Suppose first that $F_K$ is obtained from $F_K'$ via tunneling. If 
 $F_K'$ satisfies the first condition, then tunneling does not 
 interfere with the c-disk and does not introduce 
 intersections with  the spine $\Sigma_{(A,K)}$. If $F_K'$ satisfies 
 the second condition, note that $d(f,f')\leq 1$ for every $f \in 
 F_K\cap P_K$ that is essential on $P_K$ and 
 $\chi(F_K')\geq\chi(F_K)+1$. The result follows by the triangle 
 inequality.
 
 If $F_K$ is obtained from $F_K'$ via tubing with all tubes contained 
 in $Y_K'$, these tubes do not affect a c-disk for $F_K'$ contained 
 in $X_K'$ and are disjoint from any spine 
 $\Sss_{(A,K)}$. Thus if $F_K'$ satisfies the 
 first condition, so does $F_K$. If $F_K'$ satisfies the second 
 condition, the curves of $P_K\cap F_K'$ are 
 not altered by the tubing and $1-\chi(F_K)\geq 1-\chi(F_K')$  so for any 
 curve essential curve $f \in F_K\cap P_K$, $d(f,\mcA) \leq 1-\chi(F_K') \leq 1-\chi(F_K) $ as desired.
 
 \end{proof}
 
The rest of this section will be dedicated to the proof of the 
following theorem.

 \begin{thm} \label{thm:key}
     Let $A_K$ be a $K$-handlebody with $\bdd A=P$ such that 
   if $P$ is a sphere, then $P_K$ has at least six 
     punctures.
     Suppose $F_K \subset A_K$ satisfies the following 
      conditions:
      
      \begin{itemize}
\item $F_K$ has no closed components, 
\item $F_K$ is c-bicompressible and c-weakly incompressible,
\item $F_K$ has no disk components,
\item all curves of $P_K \cap F_K$ are mutually essential unless they bound punctured disks on both
surfaces. 
\end{itemize}	  
      
      Then at least one of the following holds:

	  \begin{itemize}
	      \item There is a spine $\Sigma_{(A,K)}$ entirely
	      contained in $X_K$ say and $F_K$ has a c-disk in $X_K$
	      disjoint from that spine, or
	      
	      \item $d(f, \mcA) \leq 1-\chi(F_K)$ for every $f \in
	      F_K\cap P_K$ that is essential on $P_K$ unless $f$ is
	      the boundary of a $P_K$-parallel annulus component of
	      $F_K$.
	     \end{itemize}
     \end{thm}
 
     \begin{proof}

	 If c-disks for $F_K$ were incident to two different components of
	 $F_K$, then there would be a pair of such disks on opposite sides
	 of $F_K$ with disjoint boundaries violating c-weak
	 incompressibility.  So we deduce that all c-disks for $F_K$ are
	 incident to at most one component $S_K$ of $F_K$.  $S_K$
	 cannot be an annulus, else the boundaries of c-disks in $X_K$ and 
	 $Y_K$
	 would be parallel and so could be made disjoint. In
	 particular $S_K$, and thus $F_K$, must have a strictly negative Euler characteristic.

	 Suppose $F_K$ is a counterexample to the theorem such that 
	$1-\chi(F_K)$ is minimal amongst all such counterexamples. As
	 in Proposition \ref{prop:essential} we may assume that $F_K$ has no 
	 components that are $P_K$-parallel
	 annuli or $P_K$-parallel punctured disk components. In 
	 particular this implies that all curves of $F_K \cap P_K$ are 
	 mutually essential.
	 We will prove the theorem in a sequence of lemmas. We will use
	 the following definition modeled after the definition of a 
	 strongly $\bdd$-compressible surface first introduced in 
\cite{Sc2}. 
	 
     \begin{defin}
A splitting surface $F_K \subset A_K$ splitting $A$ into submanifolds 
$X$ and $Y$ is called strongly 
$P$-compressible if there exist $P$-compressing disks $E_X 
\subset X_K$ and $E_Y\subset Y_K$ for $F_K$ such that $\bdd E_X \cap \bdd 
E_Y=\emptyset$.	 
	 
\end{defin}

\begin{lemma}\label{lem:stronglycompressible}
       The surface $F_K$ that provides a counterexample to  
Theorem \ref{thm:key}
       with maximal Euler characteristic is not strongly 
$P$-compressible.

	\end{lemma}    
     
 \begin{proof}
    
By way of contradiction suppose $E_{X} \subset X_K$ and $E_{Y} \subset
Y_K$ is a pair of disjoint $P$-compressing disks for $F_K$.  Let
$F_K^{x}, F_K^{y}$ denote the surfaces obtained from $F_K$ by
$P$-compressing $F_K$ along $E_{X}$ and $E_{Y}$ respectively, and let
$F_K^{-}$ denote the surface obtained by $P$-compressing along both
disks simultaneously.  A standard innermost disk, outermost arc
argument between $E_X$ and a c-disk for $F_K$ in $X_K$ shows that
$F_K^{x}$ has a c-disk lying in $X_K$.  Similarly, $F_K^{y}$ has a
c-disk lying in $Y_K$.  If one of $F_K^x$ or $F_K^y$ has c-disks on both
sides, say $F_K^x$, then all curves of $P_K\cap F_K^x$ must be
mutually essential unless they bound punctured disks on both surface
by Proposition \ref{prop:preserving}.  $F_K^x$ cannot be the union of 
punctured disks as it is bi-compressible so at least one component 
of $F_K^x \cap P_K$ is essential on $P_K$. As $1-\chi(F_K^x)<1-\chi(F_K)$
the surface $F_K^x$ satisfies one of the conclusions of the theorem. 
By Proposition \ref{prop:recovery} tunneling to recover $F_K$ from
$F_K^x$ preserves either of these properties so $F_K$ is not a
counterexample as we assumed.

If $F_K^-$ has any c-disk, then one of $F_K^x$ or $F_K^y$ has c-disks
on both sides as c-disks are preserved under tunneling and we are done
as above.  Suppose some curve of $F_K^-\cap P_K$ is inessential on
$P_K$ but essential on $F_K^-$.  This curve must be adjacent to the
dual arc to one of the $P$-compressing disk, say the dual arc to
$E_X$.  In this case, by an argument similar to the proof of
Proposition \ref{prop:preserving}, $F_K^y$ is c-compressible in $X_K$. 
As we saw that $F_K^y$ is c-compressible in $Y_K$, it follows that $F_K^y$ c-bicompressible, a case
we have already considered.  Thus all curves essential on $F_K^-$ are
also essential on $P_K$, therefore if $F_K^-$ has a component that is
not $P_K$-parallel, the result follows from Proposition
\ref{prop:essential}.

     We have reduced the proof to the case that $F_K^-$ is
    c-incompressible, each component of $F_K^-$ is $P_K$-parallel and
    all curves of $P_K\cap F_K^-$ are essential on $P_K$ or mutually
    inessential.  It is clear that in this case we can isotope $F_K^-$
    to be disjoint from any spine $\Sigma_{(A,K)}$.  The original
    surface $F_K$ can be recovered from $F_K^-$ by tunneling along two
    arcs on opposite sides of $F_K^-$.  The tunnels can be made
    disjoint from $\Sigma_{(A,K)}$ and thus $F_K$ can also be isotoped
    to be disjoint from $\Sigma_{(A,K)}$.  Without loss of generality
    we will assume $\Sigma_{(A,K)} \subset X_K$, thus it suffices to
    show that $F_K$ has a c-disk in $X_K$ that is disjoint from
    $\Sigma_{(A,K)}$. 
   
    Consider how $F_K^x$ can be recovered from
    $F_K^-$; the $P$-compression into $Y_K$ must be undone via a
    tunnelling along an arc $\gamma$ where the
    interior of $\gamma$ is disjoint from $F_K^-$.  Let $\gamma$
    connect components $F_K^0$ and $F_K^1$ (possibly $F_K^0=F_K^1$) of
    $F_K^-$ where $F_K^i$ is parallel to a subsurface $\tilde F_K^i
    \subset P_K$.  There are three cases to consider.
     
     First assume that $F_K^0 \neq F_K^1$ and they are nested, i.e.
     $\tilde F_K^0\subset \tilde F_K^1$.  Consider the eyeglass curve
     $e=\eta(\gamma \cup \omega)$ where $\omega\subset F_K^0$ is
     parallel to the boundary component of $F_K^0$ that is adjacent to
     $\gamma$.  Using the product structure between $F_K^0$ and
     $F_K^1$, a neighborhood of $e \times I$ contains the desired
     compressing disk for $F_K$ that is disjoint from some spine $\Sigma_{(A,K)}$.
     
     Next suppose $F_K^0\neq F_K^1$ and they are not nested.  Then
     each component of $F_K^x$ is $P_K$-parallel.  As we have already
     seen, $F_K^{x}$ has a c-disk in $X_K$.  The c-disk is either
     disjoint from some $\Sigma_{(A,K)}$, in which case we are done, or,
     via the parallelism to $P_K$, the c-disk represents a c-disk $D^*$
     for $P_K$ in $A_K$ whose boundary is disjoint from at least one
     curve in $\bdd F_K^{x}$; the curve that is in the boundary of 
     the c-compressible component of $F_K^x$. Call this particular curve $f^x$.

     If $\chi(F_K)<-1$ then $d(f^{x}, \bdd D^*) \leq 1$ so $d(f, \mcA)
	  \leq 3 \leq 1 - \chi(F_K)$. If
     $\chi(F_K)=-1$, then $F_K^x$ consists only of $P_K$-parallel
     annuli and punctured disks components. Let $N$ be the annulus 
     component of $F_K^x$ with boundary $f^x$ parallel to a 
     subannulus $\tilde{N} \subset P_K$. Then $f^{x}$ and $\bdd D^*$ 
     both lie in $\tilde{N}$ so $d(f^{x}, \bdd D^*) =0$. By Proposition
     \ref{prop:cutimpliescompressing} $d(\bdd D^*, \mcA)\leq 1$. Thus for $f$ any essential
     component of $\bdd F_K$, $d(f, \mcA) \leq d(f, f^{x}) +d(f^{x},
     \bdd D^*) +d(\bdd D^*, \mcA)\leq 1+1 = 1 - \chi(F_K)$.

     The last case to consider is the case $F_K^0=F_K^1$.  If $\gamma
     \subset \tilde F_K^0$ then $\gamma \times I$ is the desired
     compressing disk.  If $\gamma$ is disjoint from $\tilde F_K^0$,
     then each component of $F_K^x$ is $P_K$-parallel.  Proceed as in
     the previous case to show that either $F_K^x$, and thus $F_K$,
     has a c-disk disjoint from $\Sigma_{(A,K)}$ or $d(f,
     \mcA) \leq 1 - \chi(F_K)$.

 \end{proof}

\begin{lemma} \label{lem:keycompressing}
If the surface $F_K$ that provides a counterexample to Theorem
\ref{thm:key} with maximal Euler characteristic is bicompressible,
then the surfaces $F_K^X$ and $F_K^Y$ obtained from $F_K$ by maximally
compressing $F_K$ into $X_K$ and $Y_K$ respectively have cut-disks.
\end{lemma}
    
    \begin{rmk}
	Note that the hypothesis of this lemma holds when $F_K$ does
	not have any cut-disks.
    \end{rmk}
    \begin{proof}

Suppose $F_K^X$ say has no-cut disks. 	
By Corollary \ref{cor:maxincomp} the surfaces $F_K^X$ and $F_K^Y$ are 
incompressible in
$A_K$.  If some component of $F_K^X$ is not $P_K$-parallel, 
then the second conclusion of the theorem follows from Proposition
\ref{prop:essential}.  We may therefore assume that there is some
spine $\Sigma_{(A,K)}$ that is disjoint from $F_K^X$.
 
  Let $X_K^-$ and $Y_K^+$ be the two sides of $F_K^X$ and let
  $\Gamma \subset Y_K^+$ be the graph dual to the compressions we
  performed, i.e. $F_K$ can be recovered from $F_K^X$ by tunneling
  along the edges of $\Gamma$.  Note that by general position we can
  always arrange that $\Gamma$ is disjoint from any spine so in
  particular after an isotopy, $F_K\cap \Sss_{(A,K)}=\emptyset$. 
  
\medskip  
  {\bf Claim:} Recall that $S_K$ is the component of $F_K$ to which all
 c-disks for $F_K$ are incident.  To prove the lemma at hand it
 suffices to show that \begin{itemize}\item $S_K$ has a c-disk $D^*$
 on the same side of $S_K$ as the spine $\Sss_{(A,K)}$ and disjoint
 from that spine, or \item there is a compressing disk for $P_K$ whose
 boundary is disjoint from at least one curve in $\bdd S_K$, or \item
 $S_K$ is strongly $P$-compressible.  \end{itemize}

 By an innermost disk argument we may 
  isotope any c-disk for $S_K$ to be disjoint from $F_K$. 
Suppose  $S_K$ has c-disk $D^*$ on the same side of 
  $S_K$ as the spine $\Sss_{(A,K)}$ and disjoint from that
spine. 
Recall that $F_K\cap \Sss_{(A,K)}=\emptyset$ so it is sufficient to 
show that $F_K$ also has a c-disk on the same side as $\Sss_{(A,K)}$ 
but disjoint from it.

If there is a component of $F_K$ that separates $D^*$ and 
$\Sss_{(A,K)}$ than this component also separates $S_K$ and all its 
c-disks from the spine. As $S_K$ is bicompressible, we can 
always find a c-disk for $S_K$ on the same side as $\Sss_{(A,K)}$ 
and all these c-disks will be disjoint from the spine. If there is no 
such separating component, then $D^*$  is a c-disk for $F_K$ on the 
same side as $\Sss_{(A,K)}$ but disjoint from $\Sss_{(A,K)}$.

In the second case, $d(s,\mcA)\leq 1$ where $s\in \bdd S_K$ so 
$d(f,\mcA)\leq 2\leq 1-\chi(F_K)$. 

In the third case, suppose first that all components of $F_K-S_K$ 
are annuli, necessarily not $P_K$-parallel. If one of these annuli is 
$P$-compressible, $P$-compressing it results in a 
compressing disk for $P_K$ that is disjoint from $F_K$ so 
$d(f,\mcA)\leq1$. Thus we may assume that all other components of 
$F_K$ are $P$-incompressible. By an innermost disk and outermost arc 
arguments, the pair of 
strongly $P$-compressing disks 
for $S_K$ can be isotoped to be disjoint from all other components of 
$F_K$ so $F_K$ is also strongly $P$-compressing disks for $F_K$ and 
by Lemma 
\ref{lem:stronglycompressible}, $F_K$ cannot be a counterexample to 
the theorem. 

If some component of $F_K$ other than $S_K$ has a strictly negative
Euler characteristic, then $1-\chi(S_K)<1-\chi(F_K)$.  This shows that
$S_K$ is not a counterexample to the theorem, so either $d(s,\mcA)\leq
1-\chi(S_K)$ in which case $d(f,\mcA)\leq d(f,s)+d(s,\mcA)\leq
1-\chi(F_K)$ or $S_K$ has a c-disk on the same side of $S_K$ as the
spine $\Sss_{(A,K)}$ but is disjoint from it.  By repeating the
argument from the first case, we conclude that $F_K$ must also satisfy
the second conclusion of the theorem.
\medskip
  
  Note that $S_K$ is itself a c-weakly incompressible surface as 
    every c-disk for the surface $S_K$ is also a c-disk for 
    $F_K$. We will prove the lemma by showing that $S_K$ satisfies 
    one of the items in the claim above. Let $S$ split $A$ into 
    submanifolds $U$ and $V$ and $S_K^U$ be the 
    surface obtained by maximally compressing $S_K$ in $U_K$, $S_K^U$ 
    splits $A_K$ into submanifolds $U_K^-$ and $V_K^+$ and $\Gamma$ 
    is the graph dual to the compressing disk. We have already shown that 
    for some spine $\Sss_{(A,K)}$,
    $\Sss_{(A,K)}\cap F_K =\emptyset$ so in particular $\Sss_{(A,K)}\cap S_K =\emptyset$. 
    As $S^U_K$ is c-incompressible, 
    we may assume each component is $P_K$-parallel as otherwise the 
    result will follow by Proposition \ref{prop:essential}.
  We will show 
    that $S_K$ satisfies one of the conditions in the claim.
  
 If $\Sss_{(A,K)}\subset U_K^-$, then $\Sss_{(A,K)}$ is also disjoint
 from every compressing disk for $S_K$ lying in $U_K$ as it is
 disjoint from the meridional circles for the edges of $\Gamma$ and
 we have the desired result.  Thus we may assume $\Sss_{(A,K)}\subset
 V_K^+$. Let $S_K^0$ be an outermost component of $S_K^U$, i.e. a component cobounding 
a product 
 region $R_K\cong S_K^0\times I$ with $P_K$ such that $R_K \cap S^U_K=\emptyset$.

 {\bf Case 1:} Suppose for some outermost component, $R_K\subset V_K^+$. As $\Gamma
 \subset V_K^+$ and $S_K$ is connected, $S_K^0$ is the only component 
 of $S_K^U$. This implies that $\Sigma_{(A,K)}\subset R_K$ so we can 
 use the product structure to push $\Sigma_{(A,K)}$ into $U_K^-$ and 
 by the previous paragraph $S_K$ satisfies the hypothesis of the 
claim.
 
  {\bf Case 2:} Suppose the components of $S_K^U$ are nested and let
  $S_K^1$ be a second outermost component.  The region between $S_K^1$
  and the outermost components of $S_K^U$ is a product region that
  must be contained in $V_K^+$ or we can apply Case 1.  Again as $S_K$
  is connected, $V_K^+$ is also connected so in fact $V_K^+$ is a
  product region and $\Sigma_{(A,K)}\subset V_K^+$.  Again we can push
  $\Sigma_{(A,K)}$ into $U_K^-$ and complete the argument as in the
  previous case.

 {\bf Case 3:} Finally suppose that the components of $S_K^U$ are all
 outermost and all outermost regions are contained in $U_K^-$.  By
 Corollary \ref{cor:nicecompdisk}, there is a compressing disk 
 for $P_K$ that is disjoint from a complete collection of compressing
 disks for $S_K$ in $U_K$ and intersects $S_K$ only in arcs that are
 essential on $S_K^U$. Take such a disk $D$ that intersects $S_K$ 
 minimally. Consider an outermost arc of $S_K^U \cap D$
 cutting off a subdisk $D_0$ from $D$.  If $D_0\subset V_K$,
 $P$-compressing $S_K$ along $D_0$ preserves the compressing disks of
 $S_K$ lying in $U_K$ (because $D$ is disjoint from all compressing 
 disks for $S_K$ in $U_K$ by hypothesis) and also preserves the c-disks 
 lying in $V_K$ (by an innermost disk argument) so
 the result follows by induction.  If every outermost disks is
 contained in $U_K^-$, the argument of Theorem 5.4, Case 3 in
 \cite{Sc2} now carries over to show that either $S_K$ is strongly 
 $P$-compressible or there is a compressing disk for $P_K$ that is 
 disjoint from $S_K$.  We repeat the argument here for completeness.
  
 If 
 there is nesting among the arcs $D\cap S_K$ on $D$, consider a second outermost 
arc
 $\lambda_{0}$ on $D$ and let $D'$ be the disk
 this arc cuts from $D$, see Figure \ref{fig:rectangle}. If every arc of $S_K^U\cap D$ is outermost 
of $D$ let $D=D'$. Let $\Lll \subset D'$ denote the collection
 of arcs $D' \cap S_K$; one of these arcs (namely $\lambda_{0}$) will
 be on $\bdd D'$.
  Consider how a c-disk $E^*$ for $S_K$ in $V_K$ intersects $D'$.  All
 closed curves in $D' \cap E^*$ can be removed by a standard innermost
 disk argument redefining $E^*$.  Any arc in $D' \cap E^*$ must have
 its ends on $\Lll$; a standard outermost arc argument can be used to
 remove any that have both ends on the same component of $\Lll$.  If
 any component of $\Lll - \lambda_{0}$ is disjoint from all the arcs
 $D' \cap E^*$, then $S_K$ could be $P$-compressed without affecting
 $E^*$.  This reduces $1 - \chi(S_K)$ without affecting
 bicompressibility, so we would be done by induction.  Hence we
 restrict to the case in which each arc component of $\Lll -
 \lambda_{0}$ is incident to some arc components of $D' \cap E^*$.

\begin{figure}[tbh]
\centering
\includegraphics[width=0.4\textwidth]{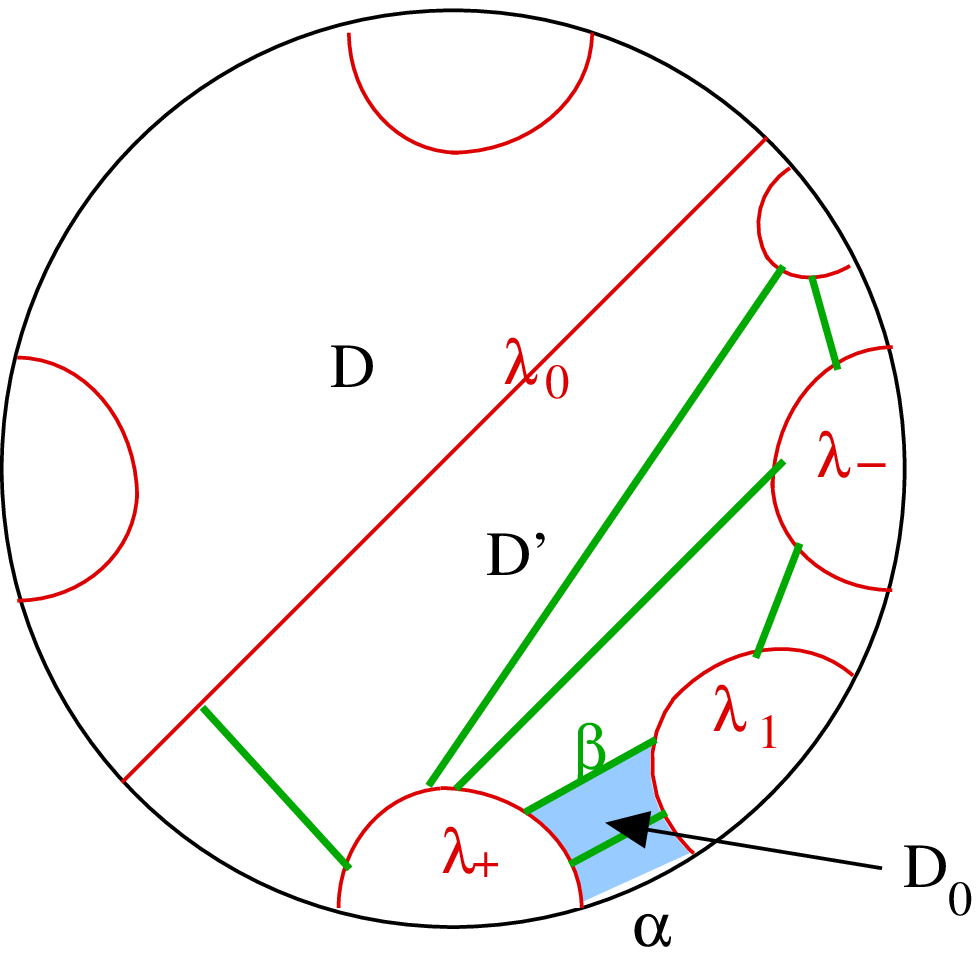}
\caption{} \label{fig:rectangle}
\end{figure}

 It follows that there is at least one component $\lambda_{1} \neq
 \lambda_{0}$ of $\Lll$ with this property: any arc of $D' \cap E^*$
 that has one end incident to $\lambda_{1}$ has its other end incident
 to one of the (at most two) neighboring components $\lambda_{\pm}$ of
 $\Lll$ along $\bdd D'$.  (Possibly one or both of $\lambda_{\pm}$ are
 $\lambda_{0}$.)  Let $\bbb$ be the outermost arc in $E^*$ among all
 arcs of $D' \cap E^*$ that are incident to the special arc
 $\lambda_{1}$.  We then know that the other end of $\bbb$ is incident
 to (say) $\lambda_{+}$ and that the disk $E_{0} \subset E^*$ cut off
 by $\bbb$ from $E^*$, although it may be incident to $D'$ in its
 interior, at least no arc of intersection $D' \cap interior(E_{0})$
 is incident to $\lambda_{1}$.  Notice that even if $E^*$ is a cut-disk, we can always choose $E_0$ so that it does not contain a
 puncture.

 Let $D_{0}$ be the rectangle in $D'$ whose sides consist of subarcs
 of $\lambda_{1}$, $\lambda_{+}$, $\bdd D'$ and all of $\bbb$. 
 Although $E^*$ may intersect this rectangle, our choice of $\bbb$ as
 outermost among arcs of $D \cap E^*$ incident to $\lambda_{1}$
 guarantees that $E_{0}$ is disjoint from the interior of $D_{0}$ and
 so is incident to it only in the arc $\bbb$.  The union of $E_{0}$ 
 and
 $D_{0}$ along $\bbb$ is a disk $D_{1} \subset V_K$ whose boundary
 consists of the arc $\aaa = P \cap \bdd D_{0}$ and an arc $\bbb'
 \subset S_K$.  The latter arc is the union of the two arcs $D_{0}
 \cap S_K$ and the arc $E_{0} \cap S_K$.  If $\bbb'$ is essential in
 $F_K$, then $D_{1}$ is a $P$-compressing disk for $S_K$ in $V_K$ that
 is disjoint from the $P$-compressing disk in $U_K$ cut off by
 $\lambda_{1}$.  So if $\bbb'$ is essential then $S_K$ is strongly
 $P$-compressible.

 Suppose finally that $\bbb'$ is inessential in $S_K$ so $\bbb'$ is
 parallel to an arc on $\bdd S_K$.  Let $D_2 \subset S_K$ be the disk
 of parallelism and consider the disk $D'=D_1 \cup D_2$.  Note that
 $\bdd D'\subset P_K$ and $D'$ can be isotoped to be disjoint from $S_K$.  Either $D'$ is
 $P_K$-parallel or is itself a compressing disk for $P_K$.  In the
 latter case $\bdd D' \in \mathcal{A}$, $d(f, \mathcal{A}) \leq 1$ for
 every $f \in \bdd S_K$ and we are done.  On the other hand if $D'$
 cobounds a ball with $P_K$, then $D_1$ and $D_2$ are parallel and so
 we can isotope $S_K$ replacing $D_2$ with $D_1$.  The result of this
 isotopy is the curves $\lambda_1$ and $\lambda_+$ are replaced by
 a single curve containing $\beta$ as a subarc lowering $|D \cap
 S_K|$.  This contradicts our original assumption that $S_K$ and $D$
 intersect minimally. We conclude that $S_K$ satisfies the second or 
 the third condition of the Claim completing the proof of Lemma 
 \ref{lem:keycompressing}.

\end{proof}

We return now to the proof of the theorem. By the above lemmas we may 
assume $F_K$ is not strongly $P$-compressible, and if it is bicompressible
both of $F_K^X$ and $F_K^Y$ have cut-disks.

     \begin{rmk}
\textnormal {Some of the argument to follow here parallels the 
argument in
Theorem 5.4 of \cite{Sc2}.  In fact it seems likely that the stronger
result proven there still holds.}
     \end{rmk}

If $F_K$ has no compressing disks on some side (and necessarily has a
cut disk), pick that side to be $X_K$.  If both sides have compressing
disks, pick $X_K$ to be the side that has a cut-disk if there is such.  Thus if $F_K$ has a cut disk, 
then it has a cut disk $D^c \subset X_K$ and if $F_K$ has a compressing disk lying
in $X_K$, it also has a compressing disk lying in $Y_K$.

Suppose $F_K$ has a cut-disk $D^c \subset X_K$. Let $\kappa$ be the component of $K-P$ that pierces through $D^c$ and
$B$ be a bridge disk of $\kappa$.  We want to consider how $F_K$
intersects $B$.  After a standard innermost disk argument, we may 
assume that the cut-disk $D^c$ intersects $B$ in a single arc $\mu$
with one endpoint lying on $\kappa$ and the other endpoint lying on a
component of $F_K \cap B$, label this component $b$ (see Figure
\ref{fig:nocompdisk}).  The curve $b$ is either a simple closed curve,
has both of its endpoints on $P_K$ or has at least 
one endpoint on $\kappa$.

\begin{figure}[tbh]
\centering
\includegraphics[width=.5\textwidth]{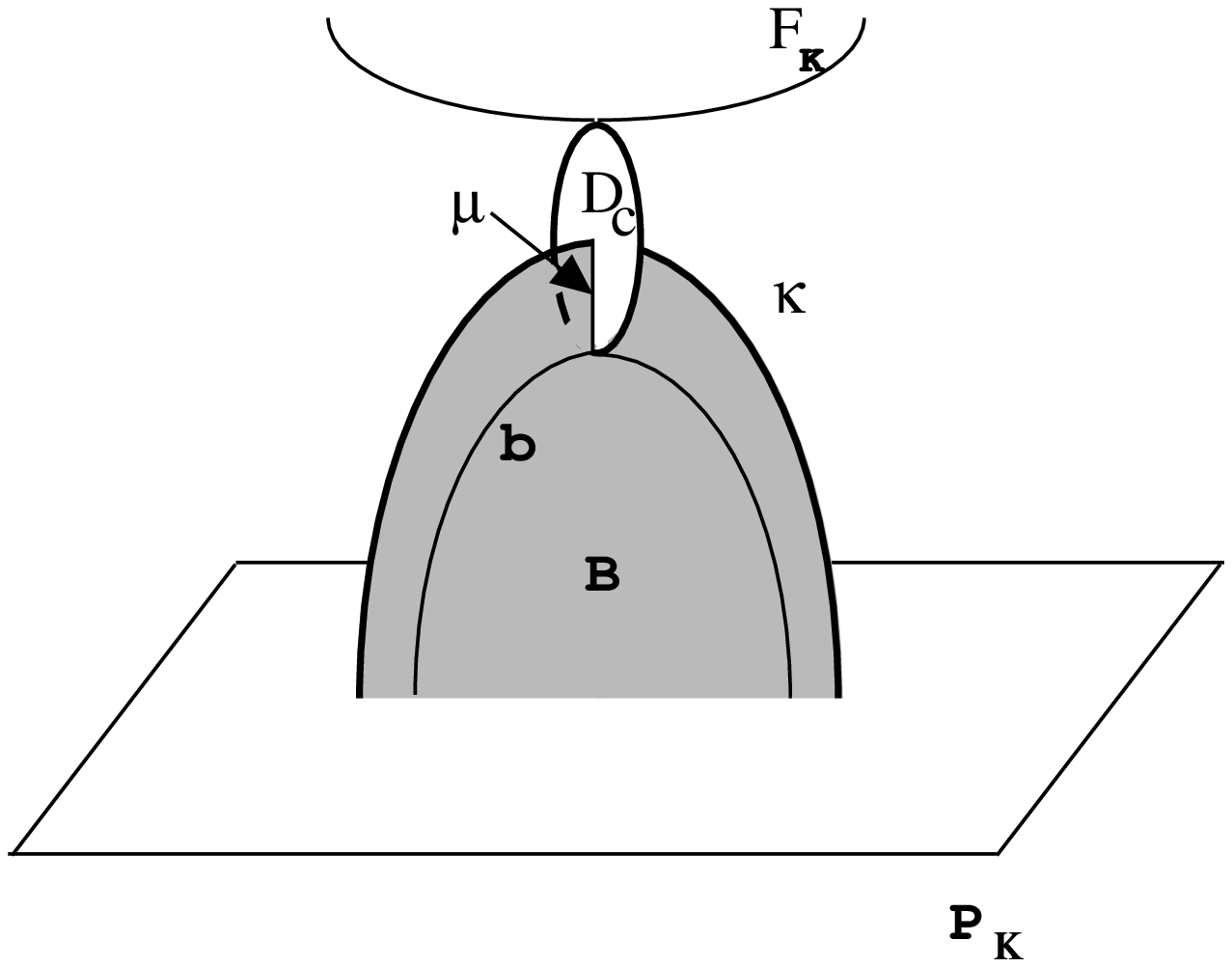}
\caption{} \label{fig:nocompdisk}
\end{figure}

Assume $|B \cap F_K|$ is minimal. We will first show that if there are
any simple closed curves of intersection, they cannot be nested on
$B$.  The argument is similar to the No Nesting Lemma in \cite{Sc3}.

Suppose such nesting occurs and let $\delta$ be a second innermost
curve cutting off a disk $D_\delta$ from $B$.  The innermost curve of
intersection contained in $D_\delta$ bound compressing disks for
$F_K$ disjoint from $D^c$ and thus must lie in $X_K$, call these 
disks $D_1,\ldots,D_n$. By our 
choice of labels this implies that $F_K$ is in fact bicompressible, 
let $E$ be a compressing disks for $F_K$ lying in $Y$. By c-weak 
incompressability of $F_K$, $E \cap D_i \neq \emptyset$. By using 
edgeslides guided by $E$ as in the proof of Lemma \ref{lem:edgeslides} 
$|B\cap F_K|$ can be reduced contradicting minimality.

We can in fact assume that there are no simple closed curves of
intersection between $F_K$ and the interior of $B$.  Suppose
$\sigma\neq b$ is an innermost simple closed curve of intersection
bounding a subdisk $D_{\sigma}\subset B$.  This disk is a compressing
disk for $F_K$ disjoint from $D^c$ so must lie in $X_K$ by c-weak
incompressibility of $F_K$.  Thus $F_K$ must also have a compressing
disk in $Y_K$.  Use this compressing disk and apply Lemma
\ref{lem:edgeslides} with the subdisk of $B$ bounded by $b$ playing the role of $T$ to isotope 
$F_K$ so as to remove all such closed
curves.

Suppose $b$ is a simple closed curve. Let $D_b\subset B$ be the disk 
$b$ bounds on $B$.  Then by the above
$D_b\cap F_K=b$ and thus $D_b$ is a compressing disk for $F_K$ lying
in $Y_K$ intersecting $D^c$ in only one point contradicting Proposition \ref{prop:nosinglepoint}. Thus we may assume $b$ is an arc.

\medskip

{\bf Case 1:} There exists a cut disk $D^c \subset
X_K$ such that the arc $b$ associated to it has both of its endpoints 
on $P_K$. 

Again let $D_b\subset B$ be the disk $b$ bounds on $B$.  By the above
discussion $D_b \cap F_K$ has no simple closed curves.  Let $\sigma$
now be an outermost in $B$ arc of intersection between $F_K$ and $B$
cutting from $B$ a subdisk $E_0$ that is a $P$-compressing disk for
$F_K$.

{\em Subcase 1A:} $b=\sigma$ and so necessarily $E_0 \subset Y_K$.  This 
in fact implies that $F_K \cap B=b$. For suppose there is an arc in 
$F_K\cap 
(B-D_b)$. An outermost such arc $\gamma$ bounds a $P$-compressing 
disk for $F_K$. If this disk is in $X_K$, then $F_K$ would be 
strongly $P$-compressible, a possibility we have already eliminated. 
If the disk is in $Y_K$, note that we can $P$-compress $F_K$ along 
this disk preserving all c-disks for $F_K$ lying in $Y_K$ and also 
preserving the disk $D^c$. The theorem then follows by Proposition 
\ref{prop:recovery}.

Consider the surface $F_K'$ obtained from $F_K$ via 
$P$-compression along $D_b$ and
the disk $D_B$ obtained by doubling $B$ along $\kappa$, a compressing 
disk for $P_K$. 
In this case $F_K'\cap D_B=\emptyset$ so we can obtain the inequality
$d(f, \mcA)\leq d(f, \bdd D_B) \leq d(f,
f')+d(f', \bdd D_B)\leq 2$ for every curve $f\in P_K\cap F_K$ as long 
as we can find at least one $f' \in
P_K\cap F_K'$ that is essential on $P_K$. 

If all curves in $P_K\cap 
F_K'$ are inessential on $P_K$, there are at most two of them.
Suppose $F_K'$ has two boundary components $f'_1$ and $f'_2$ bounding 
possibly punctured disks $D_{f'_1},D_{f'_2} \subset P_K$ and 
$F_K\cap P_K$ can be recovered by tunneling between these two curve.  As 
all curves of 
$F_K\cap P_K$ are essential on $P_K$, each of $D_{f'_1}$ and 
$D_{f'_2}$ must in fact be punctured and they cannot be nested. 
Consider the curve $f_*$ that bounds a disk on $P$ and this disk 
contains $D_{f'_1}, D_{f'_2}$ and the two points of $\kappa \cap P$, 
(see 
Figure \ref{fig:boutermosttwo}). This curve is essential on $P_K$ as 
it bounds a disk with 4 punctures on one side the other side either 
does not 
bound a disk on $P$ if $P$ is not a sphere, or contains at least two 
punctures of $P_K$ if $P$ is a sphere. As $f_*$ is disjoint from both 
the curve $F_K \cap P_K$ and from at least one curve of $\mcA$, it 
follows that the unique curve $f \in F_K\cap P_K$ satisfies the 
equality $d(F_K\cap P_K,\mcA)\leq 2\leq 1-\chi(F_K)$

\begin{figure}[tbh]
\centering
\includegraphics[width=.4\textwidth]{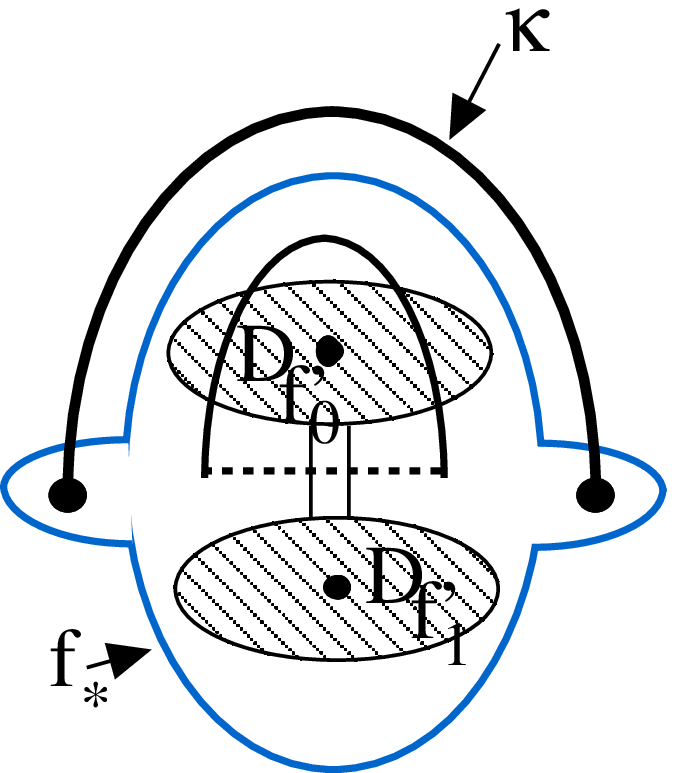}
\caption{} \label{fig:boutermosttwo}
\end{figure}

If $F_K'$ has a unique boundary curve $f'$ then $F_K$ is recovered by 
tunneling along an arc $e_0$ with both of its endpoints on $f'$. 
Therefore $F_K$ has exactly two boundary curves $f_0, f_1$ that 
cobound a 
possibly once punctured annulus on $P_K$ (see Figure 
\ref{fig:boutermost}). 

\begin{figure}[tbh]
\centering
\includegraphics[width=.6\textwidth]{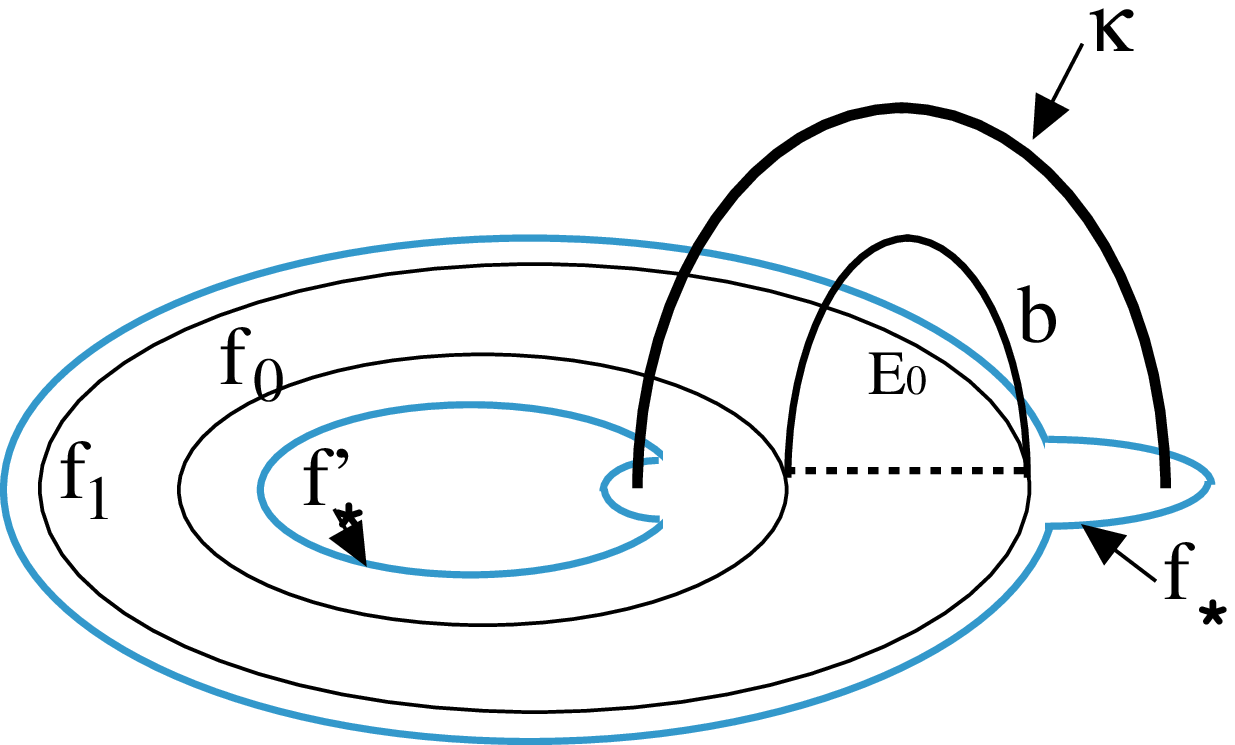}
\caption{} \label{fig:boutermost}
\end{figure}  
 
Let $f_*$ and $f'_*$ be the curves on $P_K$ that cobound once 
punctured annuli with $f_1$ and $f_0$ respectively as in Figure 
\ref{fig:boutermost}. If both $f_*$ and $f'_*$ are inessential on 
$P_K$, then $P_K$ is a sphere with at most four punctures contrary to 
the hypothesis. Thus we may assume that $f_*$ say is essential on $P_K$.  In 
this 
case $d(f_i, \mcA)\leq d(f_i, 
f_*)+d(f_*, \mcA)\leq 2\leq 1-\chi(F_K)$ for $i=1,2$ as desired.

{\em Subcase 1B:} $b \neq \sigma$ and some disk $E_0 \subset D_b$ 
bound by an outermost arc of $F_K \cap D_b$ is contained in $Y_K$. 
(It can be shown that as in subcase 1A, $F_K\cap (B-D_b)=\emptyset$ but we 
won't need this observation). $P$-compressing via $E_0$ results in a 
surface $F_K'$ with c-disks on
both sides as $E_0$ is disjoint from $D^c$.  By Proposition
\ref{prop:preserving} $F_K'$ satisfies the hypothesis and thus the 
conclusion of the theorem 
at hand and by Proposition \ref{prop:recovery} so does $F_K$ 
contradicting our assumption that $F_K$ is a counterexample. 

{\em Subcase 1C:} All outermost arcs of $F_K\cap D_b$ bound 
$P$-compressing
disks contained in $X_K$.  Consider a second outermost arc
$\lambda_{0}$ on $B$ (possibly $b$) and let $D'$ be the disk
this arc cuts from $B$.  Let $\Lll \subset D'$ denote the collection
of arcs $D' \cap F_K$; one of these arcs (namely $\lambda_{0}$) will
be on $\bdd D'$. The argument is now identical to Case 3 of Lemma 
\ref{lem:keycompressing}, and shows that $F_K$ is strongly $P$-
compressible, a possibility we have already eliminated, or 
$d(f,\mcA)\leq1$. See Figure \ref{fig:strongcomp} for 
the pair of strongly $P$-compressing disks in this case.

\begin{figure}[tbh]
\centering
\includegraphics[width=1\textwidth]{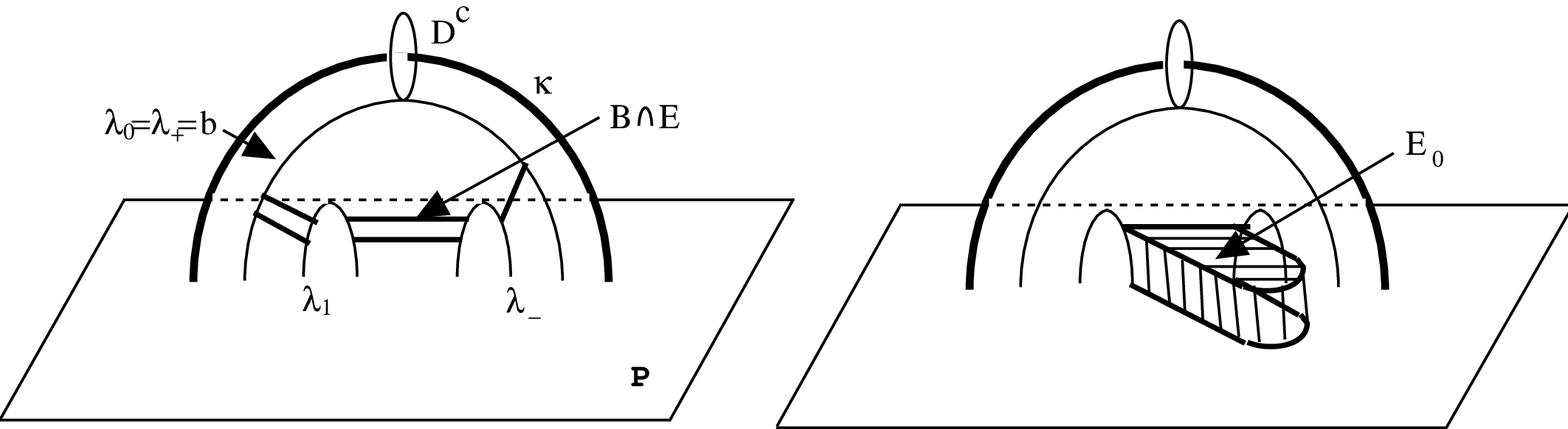}
\caption{} \label{fig:strongcomp}
\end{figure}  

\medskip

{\bf Case 2:} No cut-disk for $F_K$ has the property that the arc associated to 
it has both of its endpoints on $P_K$. 
In other words, every arc $b$ associated to a cut disk $D^c \subset 
X_K$ has at least one of its endpoints on $\kappa$. This also includes the case when $F_K$ has no cut-disks at all.

  First we will show that $F_K$ actually has compressing disks on 
  both sides. This is trivial if $F_K$ has no cut-disks so suppose $F_K$ has a cut-disk. Consider the triangle $R\subset B$ cobounded by
 $\mu, \kappa$ and $b$ (See Figure \ref{fig:assocdisk}).  If $R$ is
 disjoint from $F_K$, a neighborhood of $D^c \cup R$ contains a
 compressing disk $D$ for $F_K$, necessarily contained in $X_K$.  If 
$R \cap F_K \neq \emptyset$,
 there are only arcs of intersection as all simple closed curves have
 been removed.  An outermost on $R$ arc of intersection has both of
 its endpoint lying on $\kappa$ and doubling the subdisk of $R$ it
 cuts off results in a compressing disk $D$ for $F_K$ that also has
 to lie in $X_K$ as its boundary is disjoint from $D^c$. These two 
 types of disks will be called {\em compressing disks associated to $D^c$}.  
 As $F_K$ has a compressing disk in $Y_K$ by our initial choice of 
labeling, $F_K$ is bicompressible.

  Compress $F_K$ maximally in $X_K$ to obtain a surface $F_K^X$.  The
 original surface $F_K$ can be recovered from $F_K^X$ by tubing along 
a
 graph $\Gamma$ whose edges are the cocores of the compressing disks
 for $F_K$ on the $X_K$ side. By Corollary \ref{cor:maxincomp} 
$F_K^X$ 
 does not 
 have any compressing disks and by Lemma \ref{lem:keycompressing} it 
 has cut-disks.

 We will use $X_K^-$ and $Y_K^+$ to denote 
 the two sides of $F_K^X$ and will show that in this case $F_K^X$ 
doesn't have any cut disks lying 
in $X_K^-$. Suppose ${D'}^c \subset X_K^-$ is a cut disk for $F_K^X$ 
and $B',b'$ are respectively the disk and the arc of $F_K^X\cap B'$ 
associated to it. Note that $b'$ must have both of its endpoints on 
$P_K$ as otherwise we can construct a compressing disk associated to 
${D'}^c$ and we have shown that $F_K^X$ is incompressible. The 
original 
surface $F_K$ can be recovered from $F_K^X$ by tubing along the 
edges of a graph $\Gamma \subset Y_K^+$. This operation preserves 
the disk ${D'}^c$ and $b'$ so $F_K$ also has a cut disk whose 
associated 
arc has both of its endpoints of $P_K$ contradicting the hypothesis 
of this case. 

 \begin{figure}[tbh]
 \centering
 \includegraphics[width=1\textwidth]{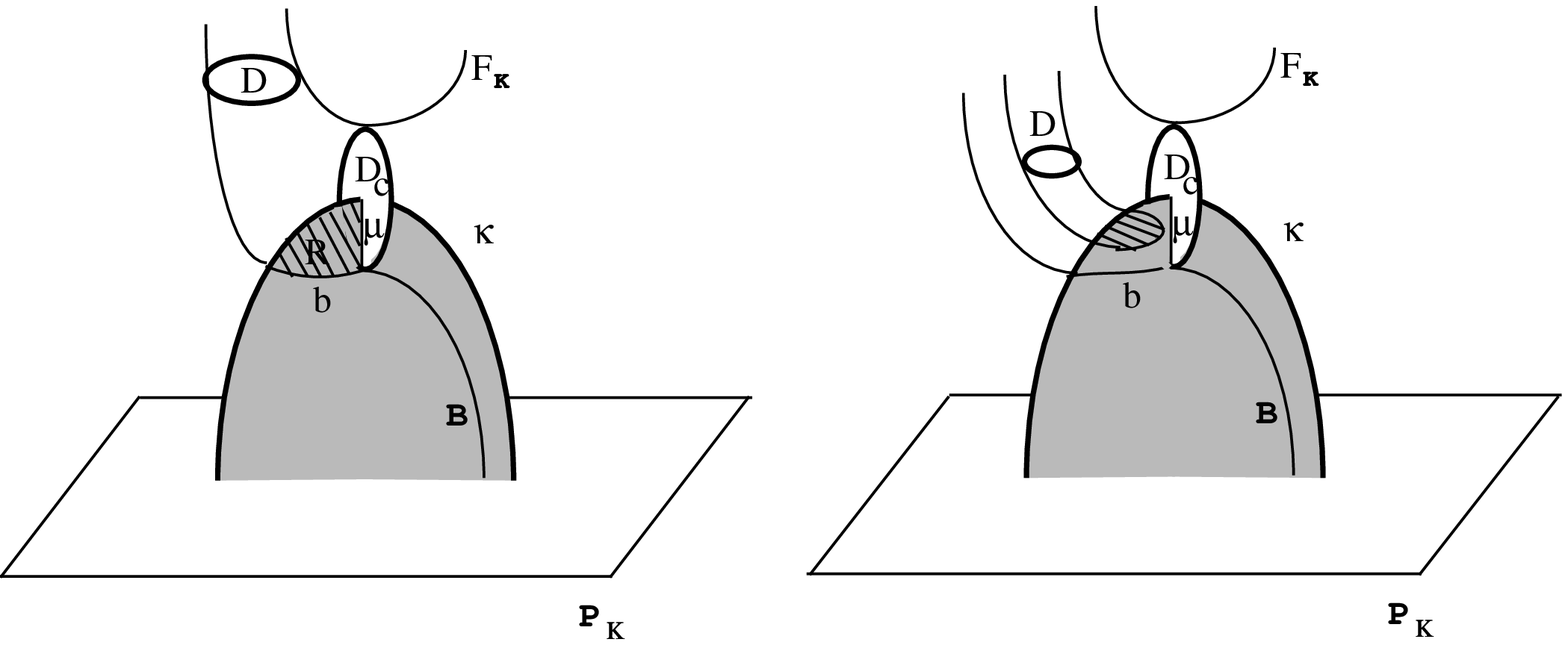}
 \caption{} \label{fig:assocdisk}
 \end{figure}

The remaining possibility is that $F_K^X$
has a cut-disk in ${D'}^c\subset Y_K^+$. Let $B'$ is its
associated bridge disk, $b'$ the arc of $F_K \cap B'$ adjacent to the
cut disk, $D_b'$ is the disk $b'$ cuts from $B'$ and $\kappa'$ the 
arc of the knot piercing ${D'}^c$.   
Assume $|F_K^X\cap B'|$ is minimal.  There cannot be any circles of
intersection for they would either be inessential on both surfaces or
give rise to compressing disks for the incompressible surface 
$F_K^X$. 
Also the arc $b'$ must have both of its endpoints on $P$, otherwise we
can construct a compressing disk for $F_K^X$ associated to ${D'}^c$ as 
in Figure \ref{fig:assocdisk}. 
Consider an outermost arc of $D_b' \cap F_K^X$ cutting from $D_b'$ a
$P$-compressing disk $E_0$, possibly $D_B'=E_0$. We now repeat an 
argument similar to the argument in Case 1 but applied to $F_K^X$. 
There are again 3 cases to consider.

{\em Subcase 2A:} 
$F_K^X \cap D_b'=b'$ so $b'$ bounds a $P$-compressing disk for 
$F_K^X$ lying 
in $X_K^-$. Let ${F_K'}^X$ be the surface obtained from $F_K^X$ after 
this $P$-compression. The argument of Subcase 1A now shows 
$d(f,\mcA)\leq 
2 \leq 2-\chi(F_K)$ for every $f\in F_K^X\cap P_K=F_K\cap P_K$.

{\em Subcase 2B:} 
Some $E_0$ lies in $Y_K^+$ (so $b'$ is not an
outermost arc). Pick a compressing disk $D$ for $F_K$ in $Y_K$ as in 
Corollary \ref{cor:nicecompdisk}. $P$-compressing $F_K$ along $E_0$ 
does not
affect c-disks lying in $Y_K^+$.  It also preserves all compressing 
disks for $F_K$ that lie in $X_K$ as it is disjoint from the graph 
$\Gamma$ and thus we are done by induction.

{\em Case 2C:} 
All outermost arcs of $F_K^X \cap B'$ lie in $X_K^-$. Consider a 
second outermost arc component of $(F_K^X) \cap B'$ and let $E_1 \in 
D_b' - 
F_K^X$ be the disk it bounds, necessarily $E_1 \subset Y^+$. By Lemma 
\ref{lem:edgeslides} we 
may assume that $\Gamma$ is disjoint from this disk. Let $E$ be a 
compressing disk for $F_K$ in $Y_K$. If $E \cap E_1 = \emptyset$ then 
$P$-compressing $F_K$ along an outermost disk component preserves the 
the compressing disk lying in $Y_K$ and of course preserves all 
c-disks lying in $X_K$ so we can finish the 
argument by induction. If there are arcs of intersection, we can 
repeat the argument of Case 1C to show that $F_K$ is strongly 
boundary 
compressible, a case we have already eliminated.

\end{proof}

     \section{Distance and intersections of Heegaard splittings}
     
For the remainder of this paper we will be considering the case of a
closed orientable irreducible 3-manifold $M$ containing a knot $K$  
with bridge surface $P$ such that $M=A\cup_P B$. In this section we 
also assume that if $P$ is a sphere then $P$ has at least 6 punctures. $Q$ will be either a second bridge surface for $K$ or 
a Heegaard surface for $M_K$. Let $X$ and $Y$ be the two components 
of $M-Q$. Thus if $Q$ is a Heegaard splitting for the knot exterior, 
then one of $X_K$ or $Y_K$ is a compression body and 
the other component is a handlebody. If $Q$ is a bridge surface, both $X_K$ and $Y_K$ 
are $K$-handlebodies.

Given a positioning of $P_K$ and $Q_K$ in $M_K$ let $Q_K^A$ and
$Q_K^B$ stand for $Q_K \cap A_K$ and $Q_K \cap B_K$ respectively. 
After removing all removable (Definition \ref{def:removable}) curves
of intersection, proceed to associate to the configuration given by
$P_K$ and $Q_K$ one or more of the following labels:

\begin{itemize}

     \item Label $A$ (resp $B$) if some component of $Q_K \cap
     P_K$ is the boundary of a compressing disk for $P_K$ lying in 
$A_K$ (resp $B_K$).
    
     \item Label $A^c$ (resp $B^c$) if some component of $Q_K \cap
	 P_K$ is the boundary of a cut disk for $P_K$ lying in $A_K$ (resp 
$B_K$).
	 
	 (Notice that this labeling is slightly different than the 
	      labeling in Section \ref{sec:essentialbound} where the 
compressing 
	      disk was required to be a subdisk of $Q_K$.)
	      
          \item $X$ (resp $Y$) if there is a compressing disk for 
$Q_K$
     lying in $X_K$ (resp $Y_K$) that is disjoint from $P_K$ and the
configuration does not already have labels $A$, $A^c$, $B$ or $B^c$.

     \item $X^c$ (resp $Y^c$) if there is a cut disk for $Q_K$ lying
     in $X_K$ (resp $Y_K$) that is disjoint from $P_K$ and the 
configuration
     does not already have labels $A$, $A^c$, $B$ or $B^c$.

     \item $x$ (resp $y$) if some spine $\Sss_{(A,K)}$ or $\Sss_{(B,K)}$ lies
     entirely in $Y_K$ (resp $X_K$) and the configuration does not 
already
     have labels $A$, $A^c$, $B$ or $B^c$.

     \end{itemize}
We will use the superscript $^*$ to denote the possible presence of
superscript $^c$, for example we will use $A^*$ if there is a label 
$A,
A^c$ or both.  

\begin{rmk} \label{rmk:inherit} If all curves of $P_K \cap Q_K$ are 
mutually
essential, then a curve is essential on $Q_K^A$ say, only if it
is essential on $Q_K$ so any c-disk for $Q_K^A$ or $Q_K^B$ that is 
disjoint from $Q_K$ is in fact
a c-disk for $Q_K$.
\end{rmk}
     
      \begin{lemma} \label{lemma:nolabel} If the configuration of $P_K$
     and $Q_K$ has no labels, then $d(K,P) \leq 2 - \chi(Q_K)$.
     \end{lemma}

     \begin{proof}
	 If $P_K \cap Q_K =\emptyset$ then $Q_K \subset A_K$ say so
	 $B_K$ is entirely contained in $X_K$ or in $Y_K$, say in
	 $Y_K$.  But $B_K$ contains all spines $\Sigma_{(B,K)}$ so
	 there will be a label $x$ contradicting the hypothesis.  Thus
	 $P_K\cap Q_K \neq \emptyset$.

	 Consider the curves $P_K \cap Q_K$ and suppose some are
	 essential in $P_K$ but inessential in $Q_K$.  An innermost
	 such curve in $Q_K$ will bound a c-disk in $A_K$ or $B_K$. 
	 Since there is no label, such curves can not exist.  In
	 particular, any intersection curve that is inessential in
	 $Q_K$ is inessential in $P_K$.  Now suppose there is a curve
	 of intersection that is inessential in $P_K$.  An innermost
	 such curve $c$ bounds a possibly punctured disk $D^* \subset
	 P_K$ that lies either in $X_K$ or in $Y_K$ but, because there
	 is no label $X^*$ or $Y^*$, this curve must be inessential in
	 $Q_K$ as well.  Let $E$ be the possibly punctured disk it
	 bounds there.  We have just seen that all intersections of
	 $E$ with $P_K$ must be inessential in both surfaces, so $c$
	 is removable and would have been removed at the onset.  We
	 conclude that all remaining curves of intersection are
	 essential in both surfaces.

 As there are no labels $X^*$ or $Y^*$, $Q_K^A$ and
 $Q_K^B$ are c-incompressible.  We conclude that both surfaces satisfy the
 hypothesis of Proposition \ref{prop:essential}.  The bound on the
 distance then follows by Corollary \ref{cor:essentialdist}.
     \end{proof}

     \begin{prop} \label{prop:AandB}
If some configuration
 is labeled $A^*$ and $B^*$ then $P_K$ is c-strongly compressible.	 
\end{prop}	 

\begin{proof}
    The labels imply the presence of c-disks for 
$P_K$ that we will denote by $D^*_A$ and $D^*_B$ such that $\bdd 
D^*_A, 
\bdd D^*_B\in Q_K \cap P_K$. As $Q_K$ is embedded, either $\bdd D^*_A=\bdd 
D^*_B$ or $\bdd D^*_A\cap\bdd 
D^*_B=\emptyset$. Thus $P_K$ is c-strongly compressible.

    \end{proof}

     \begin{lemma} \label{lem:notwellseparated}
     If $P_K\cap Q_K=\emptyset$ with say $P_K \subset X_K$ (recall 
     that $X_K$ may 
     be a handlebody, a compression body or a $K$-handlebody) and $Q_K 
\subset 
A_K$, then either every compressing disk $D$ for $Q_K$ lying in $X_K$ 
intersects $P_K$ or at least one of $P_K$ and $Q_K$ is strongly 
compressible.

     \end{lemma}

     \begin{proof}
Suppose $P_K$ and $Q_K$ are both weakly incompressible and that there is a compressing disk
for $Q_K$ lying in $X_K\cap A_K$.  As $Y_K\subset A_K$ this implies 
that $Q_K$ is bicompressible in $A_K$. As $Q_K$ is weakly 
incompressible
in $M_K$, it must be weakly incompressible in $A_K$.  Compress $Q_K$
maximally in $A_K \cap X_K$ to obtain a surface $Q_K^X$ incompressible
in $A_K$ by Corollary \ref{cor:maxincomp}.  Consider the compressing 
disks 
for $P_K$ lying in $A_K$. Each of them can be made disjoint from 
$Q_K^X$ by an innermost disk argument so the surface $P_K^A$ obtained 
by maximally compressing $P_K$ in $A_K$ is disjoint from $Q_K^X$ and 
so from $Q_K$  
(see Figure \ref{fig:notwellseparated}). As $M$ has no boundary,
$P_K^A$ is a collection of spheres and of annuli parallel to $N(K)$.  
The 
surface $P_K^A$ separates $P_K$ and $Q_K$ thus $Q_K$ is entirely 
contained in a
ball or in a ball punctured by the knot in one arc. This contradicts 
the assumption that if $M=S^3$, then $K$ is at least a three bridge 
knot. 

\begin{figure}[tbh]
\centering
\includegraphics[width=.5\textwidth]{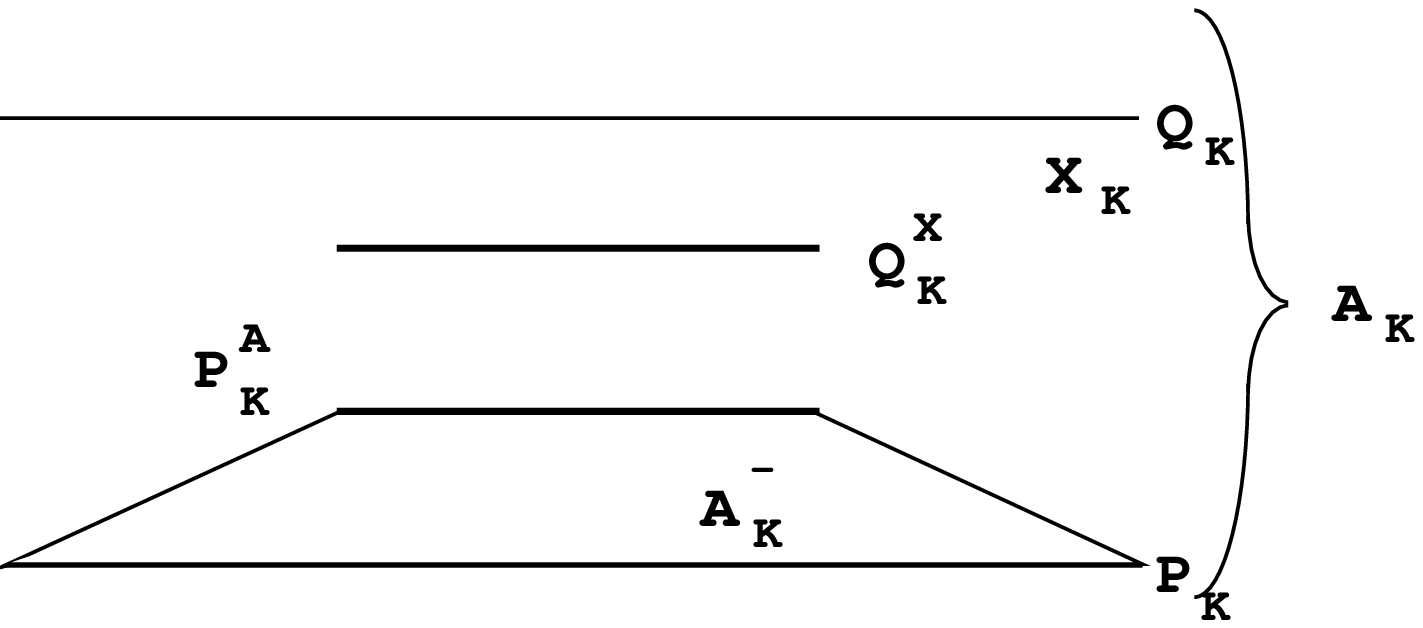}
\caption{} \label{fig:notwellseparated}
\end{figure}  
     \end{proof}

     \begin{lemma} \label{lemma:spineinY}
    If there is a spine $\Sss_{(A,K)} \subset Y_K$ (recall 
     that $Y_K$ may 
     be a handlebody, a compression body or a $K$-handlebody) then either any c-disk
     for $Q_K$ in $Y_K$ that is disjoint from $P_K$ intersects 
$\Sss_{(A,K)}$ or at least one of $P_K$ and $Q_K$ is c-strongly 
compressible.
     \end{lemma}

     \begin{proof}

      Suppose $P_K$ and $Q_K$ are both c-weakly incompressible and
    suppose $E$ is a c-disk for $Q_K$ in $Y_K$ that is disjoint from
    $P_K$ and from some spine $\Sss_{(A,K)}$.  Use the product
    structure between $P_K$ and $\Sss_{(A,K)}$ to push all of $Q_K^A$,
    as well as $E$, into $B_K$.  If $E$ was a compressing disk, this
    gives a contradiction to Lemma \ref{lem:notwellseparated} with the
    roles of $X_K$ and $Y_K$ reversed.  We want to show that even if
    the initial disk $E$ was a cut disk, after the push we can find a
    compressing disk for $Q_K$ lying in $Y_K$ that is disjoint from
    $P_K$ and contradict Lemma \ref{lem:notwellseparated}.
     
 Suppose $E$ is a cut-disk and let $\kappa\in B$ be the arc of $K-P$ that pierces $E$ and let $D
 \subset B_K$ be its bridge disk with respect to $P_K$.  Isotope $Q_K$
 and $D$ so that $|Q_K \cap D|$ is minimal and consider $b\in Q_K \cap
 D$, the arc of intersection adjacent to $E$ (this situation is
 similar to Figure \ref{fig:assocdisk}).
     
     If $b$ is a closed curve, let $ D_b$ the disk it bounds on $D$. 
     If $D \cap Q_K=b$ then $D_b$ is a compressing disk for $Q_K$ that
     intersects $E$ in exactly one point, contradicting c-weak
     incompressibility.  Let $\delta$ be an innermost curve of
     intersection between $D$ and $Q_K$ bounding a subdisk
     $D_{\delta}\subset D$.  If $D_{\delta}\subset X_K$, that would
     contradict c-weak incompressibility of $Q_K$ so $D_\delta \subset
     Y_K$ and is the desired compressing disk.

    If $b$ is not a closed curve, we can obtain a compressing disk for
    $Q_K$ much as in Figure \ref{fig:assocdisk}.  Both endpoints of
    $b$ lie on $\kappa$ as $Q_K \cap P_K =\emptyset$.  If $b$ is
    outermost, let $R$ be the disk $b$ cuts from $D$.  A neighborhood
    of $R \cup E$ consists of two compressing disks for $Q_K$ in $Y_K$
    both disjoint from $P_K$ as desired.  If $b$ is not outermost, let
    $\delta$ be an outermost arc.  Doubling the disk $D_{\delta}$ that
    $\delta$ cuts from $D$ gives a compressing disk for $Q_K$.  If
    this compressing disk is in $X_K$ that would contradict c-weak
    incompressibility of $Q_K$ thus the disk must lie in $Y_K$ as
    desired.

   \end{proof}

     Of course the symmetric statements hold if $\Sss_{(A,K)} \subset
     X_K$, $\Sss_{(B,K)} \subset Y_K$ or
$\Sss_{(B,K)} \subset X_K$.

     \begin{lemma} \label{lemma:notbothmixed}

    	   Suppose $P_K$ and $Q_K$ are both c-weakly incompressible
	    surfaces.  If there is a configuration
	   labeled both $x$ and $Y^*$ (or symmetrically $X^*$ and $y$)
	   then either $P_K$ and $Q_K$ are $K$-isotopic or $d(K,P)
	   \leq 2 - \chi(Q_K)$.
      
\end{lemma}

     \begin{proof} From the label $x$ we may assume, with no loss of
     generality, that there exists a spine $\Sss_{(A,K)} \subset Y_K$.  From
     the label $Y^*$ we know that $Q_K$ has a c-disk in $Y_K - P_K$,
     call this disk $E$.  By Lemma \ref{lemma:spineinY}, $E \cap
     \Sss_{(A,K)}\neq \emptyset$ so in particular $E \subset Y_K$.

     We first argue that we may as well assume that all components of
     $P_K \cap Q_K$ are essential in $P_K$.  For suppose not; let $c$
     be the boundary of an innermost possibly punctured disk $D^*$ in
     $P_K - Q_K$.  If $c$ were essential in $Q_K$ then $D^*$ cannot be
     in $Y_K$ (by Lemma \ref{lemma:spineinY}) and so it would have to
     lie in $X_K$.  But then $D^*$ is disjoint from $E$, contradicting
     the c-weak incompressibility of $Q_K$.  We deduce that $c$ is
     inessential in $Q_K$ bounding a possibly punctured subdisk $D'
     \subset Q_K$.  If $D'$ intersects $P_K$ in any curves that are
     essential, that would result in a label $A^*$ or $B^*$
     contradicting our labeling scheme so $c$ is removable and should
     be been removed at the onset.  Suppose now that some curve of
     intersection bounds a possibly punctured disk on $Q_K$.  By the
     above it must be essential on $P_K$ but then an innermost such
     curve would give rise to a label $A^*$ or $B^*$ contradicting the
     labeling scheme.  Thus all curves of $Q_K \cap P_K$ are mutually
     essential.

     Consider first $Q_K^B$.  It is incompressible in $B_K$ because a
     compression into $Y_K$ would violate Lemma \ref{lemma:spineinY}
     and a compression into $X_K$ would provide a c-weak compression
     of $Q_K$.  If $Q_K^B$ is not essential in $B_K$ then every
     component of $Q_K^B$ is parallel into $P_K$ so in particular
     $Q_K^B$ is disjoint from some spine $\Sigma_{(B,K)}$ and thus
     $Q_K\subset P_K\times I$.  If $Q_K$ is incompressible in
     $P_K\times I$, then it is $P_K$-parallel by Lemma
     \ref{lem:collar} as we know that $Q_K$ is not a sphere or an 
     annulus.  A compression for $Q_K$ in $P_K\times I$ would contradict Lemma
     \ref{lemma:spineinY} unless both $\Sigma_{(A,K)}$ and
     $\Sigma_{(B,K)}$ are contained in $Y_K$ and $Q_K$ has a
     compressing disk $D^X$ contained in $(P_K\times I)\cap X_K$.  In
     this case, as each component of $Q_K^B$ is $P_K$-parallel, we can
     isotope $Q_K$ to lie entirely in $A_K$ so that $P_K \subset Y_K$ 
     but then the disk $E$ provides a contradiction to  
     Lemma \ref{lem:notwellseparated}.  We conclude that
     $Q_K^B$ is essential in $B_K$ so by Proposition
     \ref{prop:essential} for each component $q$ of $ Q_K \cap P_K$
     that is not the boundary of a $P_K$-parallel annulus in $B_K$,
     the inequality $d(q, \mathcal{B}) \leq 1 - \chi(Q_K^B)$ holds. 
     Thus we can conclude that either $P_K$ and $Q_K$ are K-isotopic
     or $Q_K^B$ satisfies the hypotheses of Lemma
     \ref{lem:essentialdist}.

    By Lemma \ref{lemma:spineinY} $Q_K^A$ 
     does not have c-disks in $Y_K \cap (A_K - \Sss_{(A,K)})$ so it 
     either has 
     no c-disks in $A_K - \Sss_{(A,K)}$ at all
     or has a c-disk lying in $X_K$.  The latter would imply that
     $Q_K^A$ is actually c-bicompressible in $A_K$.  In either case we
     will show that $Q_K^A$ also satisfies the hypotheses in Lemma
     \ref{lem:essentialdist} and the conclusion of that lemma
     completes the proof.

     {\bf Case 1:} $Q_K^A$ is incompressible in $A_K - \Sss_{(A,K)}
     \cong P_K \times I$.

     By Lemma \ref{lem:collar} each component of $Q_K^A$ must be
     $P_K$-parallel.  The c-disk $E$ of $Q_K^A$ in $Y_K- P_K$ can be
     extended via this parallelism to give a c-disk for $P_K$ that is
     disjoint from all $q \in
     Q_K \cap P_K$.  Hence $d(q, \mathcal{A}) \leq 2\leq 1 -
     \chi(Q_K^A)$ as long as $Q_K^A$ is not a collection of
     $P_K$-parallel annuli.  If that is the case, then $d(\bdd E,
     q_0)$=0 for at least one $q_0 \in (P_K \cap Q_K)$ so $d(q_0,
     \mathcal{A})\leq 1 \leq 1-\chi(Q_K^A)$ as desired.

     {\bf Case 2:} $Q_K^A$ is c-bicompressible in $A_K$.  Every c-disk
     for $Q_K$ in $Y_K$ intersects $\Sigma_{(A,K)}$, so we can deduce
     the desired distance bound by Theorem \ref{thm:key}.

     \end{proof}

     \begin{lemma} \label{lemma:notboth}

	 Suppose $P_K$ and $Q_K$ are both c-weakly incompressible
	  surfaces.  If there is a
	 configuration labeled both $X^*$ and $Y^*$ then either $P_K$
	 and $Q_K$ are K-isotopic or $d(K,P) \leq 2 - \chi(Q_K)$.

	 \end{lemma}
	 
	 \begin{proof}
     Since $Q_K$ is c-weakly incompressible, any pair of c-disks, one
     in $X_K$ and one in $Y_K$, must intersect on their boundaries and
     so cannot be separated by $P_K$.  It follows that if both labels
     $X^*$ and $Y^*$ appear, the boundaries of the associated c-disks
     lie on one of $Q_K^A$ or $Q_K^B$, say, $Q_K^A$.

	  Again we may as well assume that all components of $P_K \cap
	  Q_K$ are essential in $P_K$.  For suppose not; let $c$ be
	  the boundary of an innermost possibly punctured disk $D^*$
	  in $P_K - Q_K$.  If $c$ were essential in $Q_K$ then a
	  c-disk in $B_K$ parallel to $D$ would be a c-disk for
	  $Q_K^{B}$.  From this contradiction we deduce that $c$ is
	  inessential in $Q_K$ and proceed as in the proof of Lemma
	  \ref{lemma:notbothmixed}.  As no labels $A^*$ or $B^*$
	  appear, all curves are also essential on $Q_K$.

If $Q_K^A$ or $Q_K^B$ could be made disjoint from some spine
$\Sss_{(A,K)}$ or $\Sss_{(B,K)}$, then the result would follow by
Lemma \ref{lemma:notbothmixed} so we can assume that is not the case. 
In particular $Q_K^B$ is essential and so via Proposition
\ref{prop:essential} it satisfies the hypothesis of Lemma
\ref{lem:essentialdist}.  The surface $Q_K^A$ is c-bicompressible,
c-weakly incompressible and there is no spine $\Sigma_{(A,K)}$
disjoint from $Q_K^A$.  By Theorem \ref{thm:key}, $Q_K^A$ also
satisfies the hypothesis of Lemma \ref{lem:essentialdist} so we have
the desired distance bound.
	  \end{proof}

    	 \begin{lemma} \label{lemma:notbothsmall}
 Suppose $P_K$ and $Q_K$ are both c-weakly incompressible
surfaces. If there is a configuration labeled both $x$ and $y$, then
       either $P_K$ and $Q_K$ are $K$-isotopic or $d(K,P) \leq 2 -
       \chi(Q_K)$.
      \end{lemma}

     \begin{proof}
     As usual, we can assume that all curves in $P_K \cap Q_K$ are
     essential in both surfaces.  Indeed, if there is a curve of
     intersection that is inessential in $P_K$ then an innermost one
     either is inessential also in $Q_K$, and can be removed as
     described above, or is essential in $Q_K$ and so would give a
     rise to a label $X^*$ or $Y^*$, a case done in Lemma
     \ref{lemma:notbothmixed}.  In fact we may assume that $Q_K^A$ or
     $Q_K^B$ are incompressible and c-incompressible as otherwise the
     result would follow by Lemma \ref{lemma:notbothmixed}.  As no
     labels $A^*$ or $B^*$ appear, we can again assume that all 
     curves $P_K \cap Q_K$
     are also essential on $Q_K$.

     Both $X_K$ and $Y_K$ contain entire spines of $A_K$ or $B_K$,
     though since we are not dealing with fixed spines the labels
     could arise if there are two distinct spines of $A_K$, say, one
     in $X_K$ and one in $Y_K$.  Indeed that is the case to focus on,
     since if spines $\Sss_{(A,K)} \subset X_K$ and $\Sss_{(B,K)}
     \subset Y_K$ then $Q_K$ is an incompressible surface in $P_K
     \times I$ so by Lemma \ref{lem:collar} $Q_K$ is K-isotopic to
     $P_K$.

     So suppose that $\Sss_{(A,K)} \subset Y_K$ and there is another 
     spine
     $\Sss'_{(A,K)} \subset X_K$.  $Q_K^A$ is incompressible in $A_K$
     so it is certainly incompressible in the product $A_K -
     \Sss_{(A,K)}$ and so every component of $Q_K^A$ is parallel in
     $A_K - \Sss_{(A,K)}$ to a subsurface of $P_K$.  Similarly every
     component of $Q_K^A$ is parallel in $A_K - \Sss'_{(A,K)}$ to a
     subsurface of $P_K$.

    Let $Q_{0}$ be a component of $Q_K^A$ that lies between
    $\Sss_{(A,K)}$ and $\Sss'_{(A,K)}$.  This implies that $Q_{0}$ is
    parallel into $P_K$ on both its sides, i.e. that $A_K \cong Q_{0}
    \times I$.

    As $K$ is not a 2-bridge knot, then either
    $\chi(P_K) <-2$ (so in particular $\chi(Q_0)< -1$) or $P_K$ a
    twice punctured torus.  We will show that in either case
    $d(\mathcal{A}, q) \leq 2 \leq 1-\chi(Q_K^A)$.

If $P_K$ is a twice punctured torus, then 
$Q_0$ is a once punctured annulus so has Euler characteristic -1 and 
thus $\chi(Q_K^A)<0$. Note that $d(\bdd Q_0, \mathcal{A})\leq 2$ (see 
Figure \ref{fig:torusdist}) and thus $d(\mathcal{A}, q) \leq 2
     \leq 1-\chi(Q_K^A)$.

\begin{figure}[tbh]
\centering
\includegraphics[width=.5\textwidth]{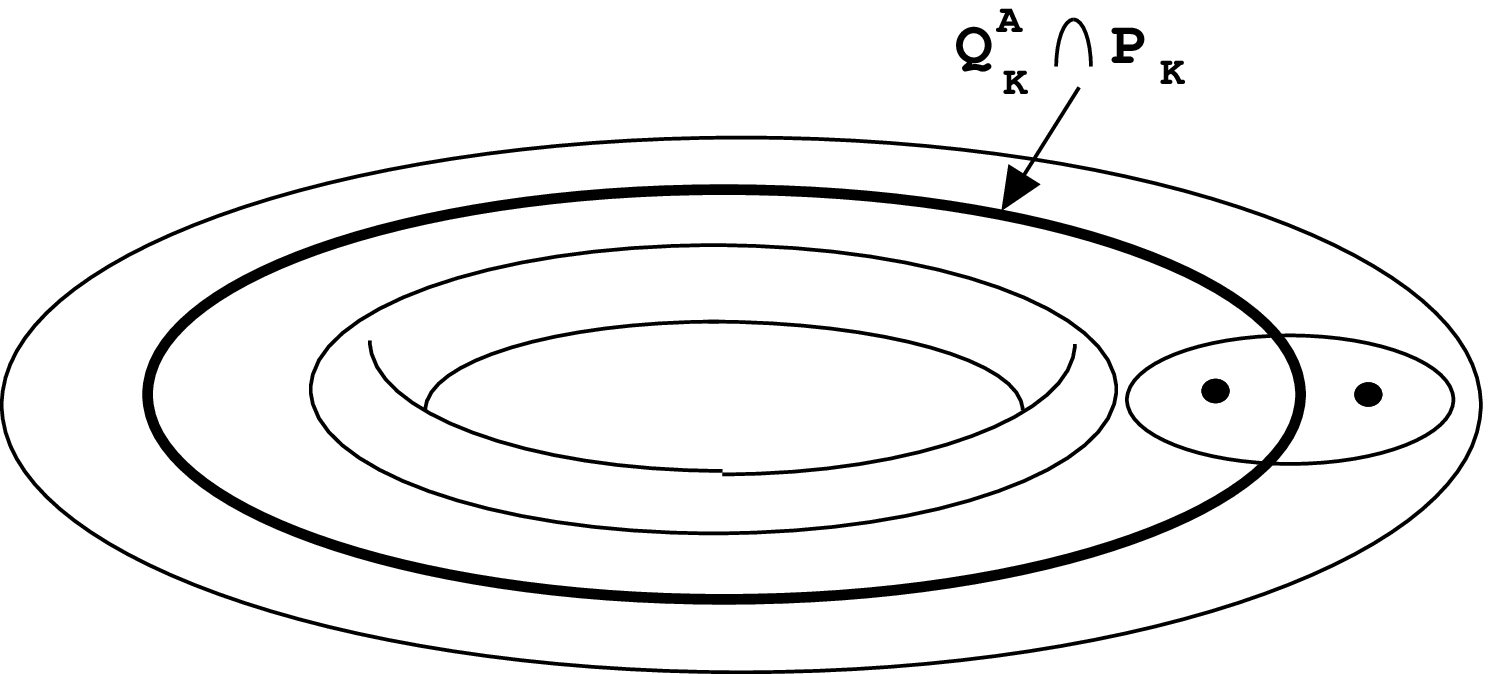}
\caption{} \label{fig:torusdist}
\end{figure}

     If $\chi(P_K) <-2$ let $\aaa$ be an essential arc in $Q_{0}$ with 
     endpoints on $P_K\cap Q_K$. 
    Then $\aaa \times I \subset Q_{0} \times I \cong A_K$ is a
    meridian disk $D$ for $A_K$ that intersects $Q_{0}$ precisely in
    $\aaa$.  $P$-compressing $Q_{0}$ along one of the two disk
    components of $D - \aaa$ produces at most two surfaces at least
    one of which, $Q_1$ say, has a strictly negative Euler
    characteristic.  In particular it is not a disk, punctured disk or
    an annulus.  Every component of $\bdd Q_1$ is essential on $P_K$
    and disjoint from both $D$ and $Q_0\cap P_K$.  We can conclude
    that for every curve $q \in P_K\cap Q_K$, $d(\mathcal{A}, q) \leq d(\mathcal{A}, \bdd Q_0)+d(\bdd
    Q_0,q)\leq 2 \leq 1-\chi(Q_K^A)$.  Thus $Q_K^A$ always satisfied
    the hypothesis of Lemma \ref{lem:essentialdist}.

     Now consider $Q_K^B$.  If it is essential, then by Proposition
     \ref{prop:essential} $Q_K^{B}$ also satisfies the hypothesis of
     Lemma \ref{lem:essentialdist} and we are done by that lemma.  If
     $Q_K^B$ has c-disks in $B_K$, we have labels $X^*$ and $y$ (or
     $x$ and $Y^*$) and we are done via Lemma
     \ref{lemma:notbothmixed}.  Finally, if each component of
     $Q_K^{B}$ is parallel to a subsurface of $P_K$, then $Q_K$ is
     disjoint from a spine $\Sigma_{(B,K)}$ as well, a case we have
     already considered.
     \end{proof}

\section{How labels change under isotopy}

Suppose $P$ and $Q$ are as defined in the previous section and 
continue to assume that if $P$ is a sphere, then $P_K$ has at least 6 punctures.
Consider how configurations and their labels change as $P_K$ say is 
isotoped while keeping $Q_K$ fixed.  Clearly if there are no 
tangencies of 
$P_K$ and $Q_K$ during the
isotopy then the curves $P_K \cap Q_K$ change only by isotopies and 
there is no change in labels.  Similarly, if there is an index
$0$ tangency, $P_K \cap Q_K$ changes only by the addition or deletion 
of a
removable curve.  Since all such curves are removed before labels are 
defined, again there is no affect on the labeling. There are two 
cases 
to consider; $P_K$ passing through a saddle tangency for $Q_K$ and 
$P_K$ passing through a puncture of $Q_K$. Consider first what can 
happen to the labeling when passing through a 
saddle tangency of $P_K$ with $Q_K$.

\begin{figure}[tbh]
\centering
\includegraphics[width=0.7\textwidth]{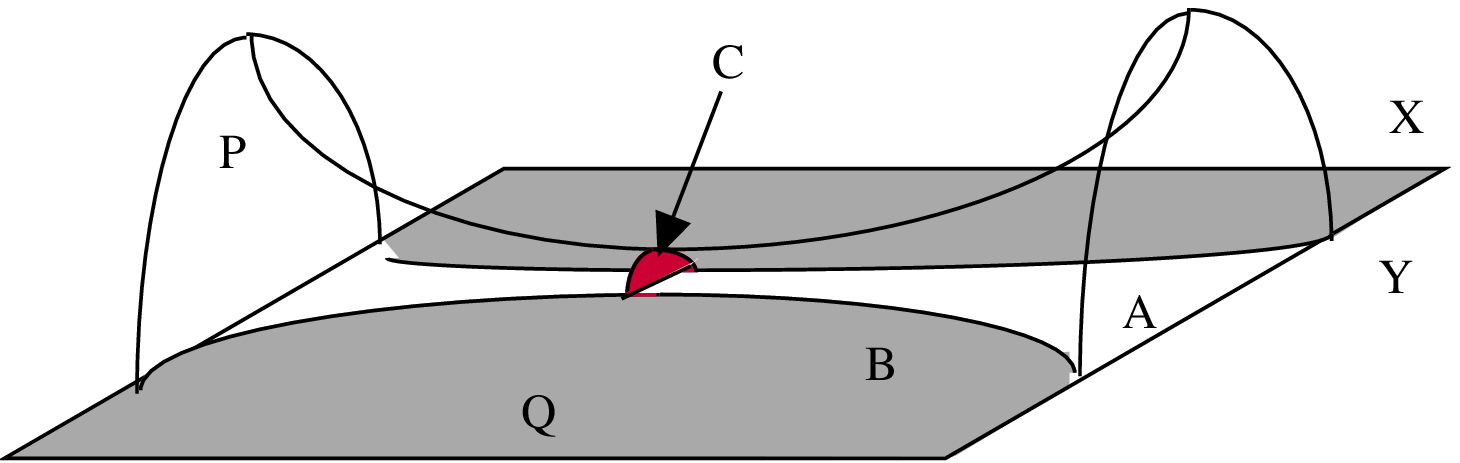}
\caption{} \label{fig:saddle}
\end{figure}

\begin{lemma} \label{lemma:saddle}
 Suppose $P_K$ and $Q_K$ are c-weakly incompressible surfaces 
  and $P_K$ is isotoped to pass through a single saddle
 tangency for $Q_K$. Suppose farther that the bigon $C$ defining the saddle tangency (see Figure 
 \ref{fig:saddle}) lies in $X_K \cap A_K$.  Then
 
 \begin{itemize}
 \item  No label $x$ or $X^*$ is removed.
 
 \item  No label $y$ or $Y^*$ is created.
 
 \item If there is no label $x$ or $X^*$ before
 the move, but one is created after and if there is a label $y$ or 
$Y^*$
 before the move and none after, then either $P_K$ and $Q_K$ are 
isotopic
 or $d(K,P) \leq 2 - \chi(Q_K)$.
 
 \end{itemize}
 
     \end{lemma}
    
\begin{proof} 
    Much of the argument here parallels the argument in the proof of 
    Lemma 4.1 in \cite{ST}. The main difference is in Claim 2.
    
    We first show that no label $x$ or $X^*$ is removed.  If 
there is a c-disk for $Q_K$ in $X_K\cap A_K$, a standard 
innermost disk, outermost arc argument on its intersection with $C$ 
shows that there is a c-disk for $Q_K$ in $X_K\cap A_K$ that is disjoint 
from $C$.  The saddle move has no effect on such a disk. It is clear 
that the move doesn't have an effect on a c-disk for $Q_K$ lying in 
$X_K \cap B_K$ so a label $X^*$ will not be removed.  If there is a spine of 
$(A,K)$ or $(B,K)$ lying entirely in $Y_K$ then that spine, too, in 
unaffected by the saddle move.
    
Dually, no label $y$ or $Y^*$ is created: the inverse saddle move, 
restoring the original configuration, is via a bigon that lies in 
$B_K 
\cap Y_K$.
    
To prove the third item position $Q_K$ so that it is exactly tangent 
to $P_K$ at the
saddle.  A bicollar of $Q_K$ then has ends that correspond to the
position of $Q_K$ just before the move and just after.  Let $Q_K^{a}$
denote $Q_K \cap A_K$ after the move and $Q_K^{b}$ denote $Q_K \cap 
B_K$ before
the move.  The bicollar description shows that $Q_K^{a}$ and $Q_K^{b}$
have disjoint boundaries in $P_K$.  Moreover the complement of 
$Q_K^{a}
\cup Q_K^{b}$ in $Q_K$ is a regular neighborhood of the singular 
component
of $P_K \cap Q_K$, with Euler characteristic $-1$.  It follows that
$\chi(Q_K^{a}) + \chi(Q_K^{b}) = \chi(Q_K) + 1$.

With $Q_K$ positioned as described, tangent to $P_K$ at the saddle 
point 
but otherwise in general position, consider the closed (non-singular) 
curves of intersection.  

{\bf Claim 1:} It suffices to consider the case in which all 
non-singular 
curves of intersection are essential in $P_K$.  

To prove the claim, suppose a non-singular curve is inessential and
consider an innermost one.  Assume first that the possibly punctured 
disk 
$D^*$ that it
bounds in $P_K$ does not contain the singular curve $s$ (i.e. the
component of $P_K \cap Q_K$, homeomorphic to a figure 8, that 
contains the
saddle point).  If $\bdd D^*$ is essential in $Q_K$, then it would 
give
rise to a label $X^*$ or a label $Y^*$ that persists from before the 
move
until after the move, contradicting the hypothesis.  Suppose on the
other hand that $\bdd D^*$ is inessential in $Q_K$ and so bounds a 
possibly punctured disk $E^*
\subset Q_K$.  All curves of intersection of $E^*$ with $P_K$ must be
inessential in $P_K$, since there is no label $A^*$ or $B^*$.  It 
follows
that $\bdd D^* = \bdd E^*$ is a removable component of intersection 
so the
disk swap that replaces $E^*$ with a copy of $D^*$, removing the 
curve of
intersection (and perhaps more such curves) has no effect on the
labeling of the configuration before or after the isotopy. So the
original hypotheses are still satisfied for this new configuration of 
$P_K$ and $Q_K$.  

Suppose, on the other hand, that an innermost non-singular inessential
curve in $P_K$ bounds a possibly punctured disk $D^*$ containing the 
singular component $s$.
When the saddle is pushed through, the number of components in $s$
switches from one $s_{0}$ to two $s_{\pm}$ or vice versa.  
All three curves are inessential in $P_K$ since they lie in the 
punctured disk $D^*$. Two of them actually bound possibly punctured 
subdisks of $D^*$ whose 
interiors are disjoint from $Q_K$. Neither of these curves can be 
essential on $Q_K$ otherwise they 
determine a label $X^*$ or $Y^*$ that persist throughout the isotopy. At least one of these curves must 
bound a nonpunctured disk on $P_K$ (as $D^*$ has at most one 
puncture) and 
thus it also bounds a nonpunctured disk on $Q_K$. We conclude 
that at least two of the curves are inessential on $Q_K$ and at least 
one of them bounds a disk on $Q_K$. As the three curves cobound a pair 
of pants of $Q_K$ the third curve is also 
inessential on $Q_K$.  This implies that all the curves are 
removable so passing through the singularity has no effect on the 
labeling. This proves the claim.

{\bf Claim 2:} We may assume that if any of the curves $s_{0}, s_{\pm}$ 
are inessential in $P_K$ they bound punctured disks on both surfaces.  
    
The case in which all three curves are inessential in $P_K$ is covered
in the proof of Claim 1.  If two are inessential in $P_K$ and at 
least 
one of them bounds a disk with no punctures then the third curve is 
also inessential. Thus if exactly two curves are inessential on 
$P_K$, they 
both bound punctured disks on $P_K$ and as no capital labels are 
preserved during the tangency move, they also bound punctured disks 
on $Q_K$ which are parallel into $P_K$. 

We are left to consider the case in which exactly one of $s_{0}, s_{\pm}$ is inessential
in $P_K$, bounds a disk there and, following Claim 1, the disk it 
bounds in $P_K$ is disjoint
from $Q_K$.  If the curve were essential in $Q_K$ then there would 
have to
be a label $X$ or $Y$ that occurs both before and after the saddle
move, a contradiction.  If the curve is inessential in $Q_K$ then it 
is
removable.  If this removable curve is $s_{\pm}$ then passing through
the saddle can have no effect on the labeling. If this removable
curve is $s_{0}$ then the curves $s_{\pm}$ are parallel in both $P_K$
and $Q_K$.  
In the latter case, passing through the saddle has the same
effect on the labeling as passing an annulus component of $P_K-Q_K$ 
across
a parallel annulus component $Q_K^{0}$ of $Q_K^{A}$.  This move can 
have
no effect on labels $x$ or $y$.  A meridian, possibly punctured disk 
$E^*$ for $Y_K$ that is
disjoint from $P_K$ would persist after this move, unless $\bdd E^*$ 
is in
fact the core curve of the annulus $Q_K^{0}$.  But then the union of 
$E^*$
and half of $Q_K^{0}$ would be a possibly punctured meridian disk of 
$A_K$ bounded by a
component of $\bdd Q_K^{0} \subset P_K$.  In other words, there would 
have
to have been a label $A^*$ before the move, a final contradiction
establishing Claim 2.

\medskip

Claims 1 and 2, together with the fact that neither labels $A^*$ nor 
$B^*$ appear, reduce us to the case in which all curves of 
intersection are
essential in both surfaces both before and after the saddle move 
except perhaps some curves which bounds punctured disks on $Q_K$ and 
on $P_K$. Let $\tilde{Q_K^a}$ and $\tilde{Q_K^b}$ be the surfaces 
left over after 
deleting from $Q_K^a$ and $Q_K^b$ any $P_K$-parallel punctured disks. 
As $Q_K^a$ and $Q_K^b$ cannot be made disjoint from any 
spine $\Sigma_{(A,K)}$ or $\Sigma_{(B,K)}$, $\tilde{Q_K^a}$ and 
$\tilde{Q_K^b}$ are not empty 
and, as we are removing only punctured disks, $\chi 
(Q_K^a)=\chi(\tilde{Q_K^a})$ 
and $\chi (Q_K^b)=\chi(\tilde{Q_K^b})$.

Note then that $\tilde{Q_K^a}$ and $\tilde{Q_K^b}$ are 
c-incompressible in 
$A_K$ and $B_K$
respectively.  For example, if the latter has a c-disk in $B_K$, then 
so 
does $Q_K^a$. Since no label $X^*$ exists before the
move, the
c-disk must be in $Y_K$ and such a c-compression would persist after the move and so then would the
label $Y^*$.  Similarly neither $\tilde{Q_K^a}$ nor 
$\tilde{Q_K^b}$ can
consist only of $P_K$ parallel components.  For example, if
all components of $\tilde{Q_K^b}$ are parallel into $P_K$ then 
$Q_K^{b}$ is also
disjoint from some spine of $B_K$ and such a spine will be unaffected 
by
the move, resulting in the same label ($x$ or $y$) arising before and
after the move.  We deduce that $\tilde{Q_K^a}$ and $\tilde{Q_K^b}$ 
are essential
surfaces in $A_K$ and $B_K$ respectively.  

Now apply Proposition \ref{prop:essential} to both sides: Let $q_{a}$
(resp $q_{b}$) be a boundary component of an essential component of
$\tilde{Q_K^a}$ (resp $\tilde{Q_K^b}$).  Then $$d(K,P) = 
d(\mathcal{A}, \mathcal{B}) \leq
d(q_{a}, \mathcal{A}) + d(q_{a}, q_{b}) + d(q_{b}, \mathcal{B}) $$ 
$$\leq 3 - \chi(\tilde{Q_K^a}) - \chi(\tilde{Q_K^b}) = \leq 3 - 
\chi(Q_K^a) - \chi(Q_K^b)=2 - \chi(Q_K)$$ as required.
    \end{proof}
 
\begin{figure}[tbh]
    \centering
    \includegraphics[width=0.7\textwidth]{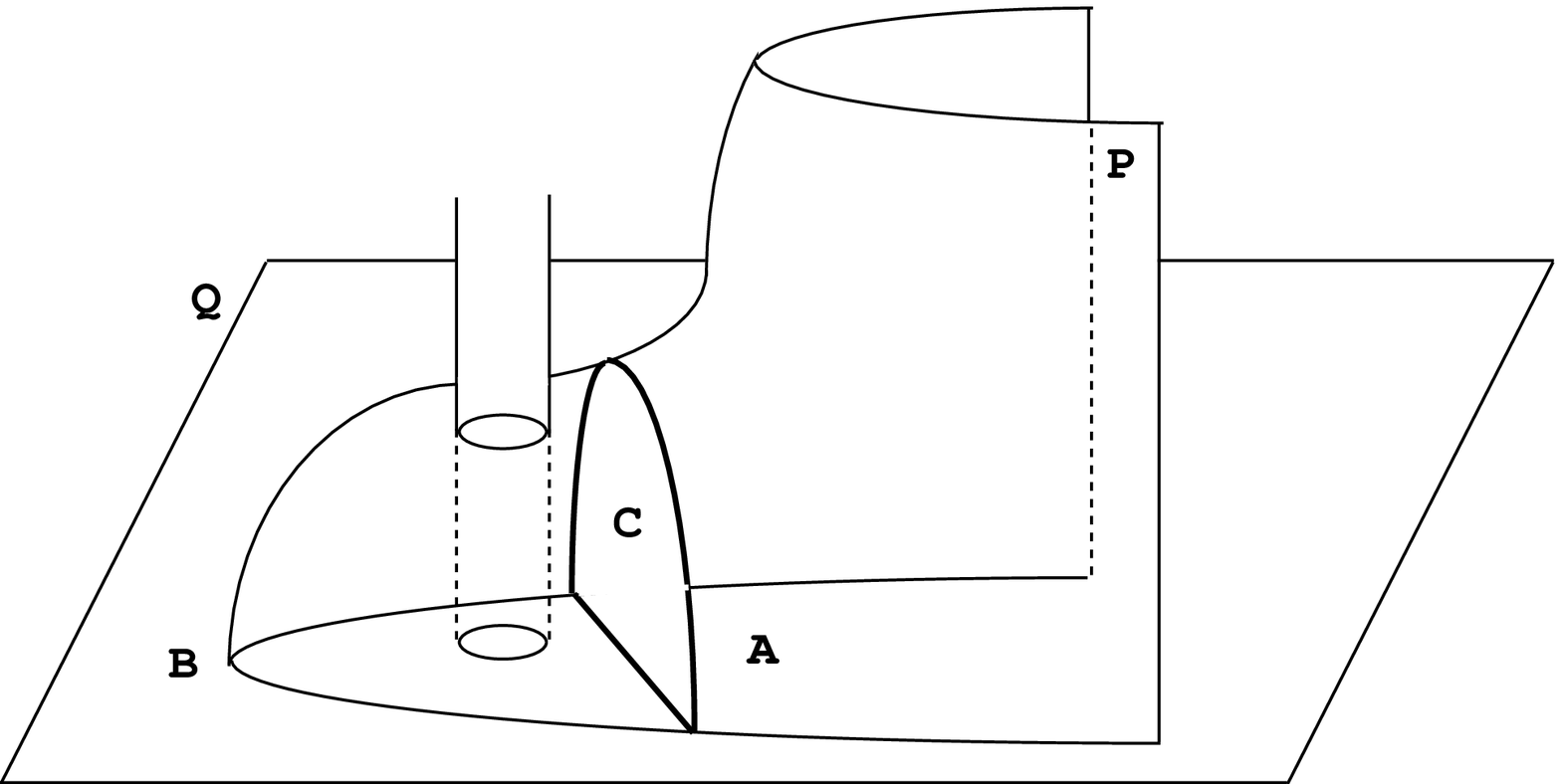}
    \caption{} \label{fig:puncture}
    \end{figure}

   It remains to consider the case when $P_K$ passes through a 
   puncture of $Q_K$ as in Figure \ref{fig:puncture}. This puncture 
   defines a bigon $C$ very similar to the tangency bigon in the 
   previous lemma: let 
   $Q_K^a$ and $Q_K^b$ be as before, then $Q_K- (Q_K^a \cup Q_K^b)$ 
is a punctured 
   annulus. The knot strand that pierces it is parallel to this 
annulus, 
   let $C$ be the double of the parallelism rectangle so that $C 
\subset 
   X_K \cap A_K$.   
    
    \begin{lemma} \label{lem:puncture}
	Suppose $P_K$ and $Q_K$ are c-weakly incompressible bridge surfaces 
	 for a knot $K$ and $P_K$ is isotoped to pass through a single 
	 puncture for $Q_K$. Suppose farther that the bigon $C$ defined by the puncture (see Figure 
	\ref{fig:puncture}) lies in $X_K\cap A_K$.  Then
	
	\begin{itemize}
	\item  No label $x$ or $X^*$ is removed.
	
	\item  No label $y$ or $Y^*$ is created.
	
	\item Suppose that, among the labels both before and after the move,
	neither $A^*$ nor $B^*$ occur.  If there is no label $x$ or $X^*$ 
before
	the move, but one is created after and if there is a label $y$ or 
$Y^*$
	before the move and none after, then either $P_K$ and $Q_K$ are 
$K$-isotopic
	or $d(P,K) \leq 2 - \chi(Q_K)$.
	
	\end{itemize}
	
	    \end{lemma}

   \begin{proof}
   The proof is very similar to the proof of the previous lemma. 
   It is clear that if there is a c-disk for $X_K$ that lies in 
$A_K$, there 
   is a c-disk that is disjoint 
   from $C$ and thus the label survives the move. If there is a spine 
of 
   $A_K$ or $B_K$ lying entirely in $Y_K$ then that spine, too, is 
   unaffected by the saddle move. The proof of the third item is 
   identical to the proof in the above lemma in the case when at 
least 
   one of the curves $s_{0}, s_{\pm}$ bounds a punctured disk on 
$Q_K$.
       
   \end{proof}

We will use $\xX$ (resp $\yY$) to denote any subset of the labels 
$x, X, X^c$ (resp $y, Y, Y^c$). The results of the last two sections 
then can be summarized as follows:

\begin{cor} \label{corollary:notboth} If two configurations are
related by a single saddle move or going through a puncture and the 
union of all labels for both
configurations contains both $\xX$ and $\yY$ then either $P_K$ and 
$Q_K$ are 
K-isotopic or $d(K,P_K) \leq 2 -
\chi(Q_K)$
    \end{cor}
    
\begin{proof} With no loss of generality, the move is as described in Lemma
\ref{lemma:saddle} or Lemma \ref{lem:puncture}.  These lemmas shows 
that 
either we have the desired bound, or
there is a single configuration for which both $\xX$ and $\yY$ 
appear.  The result then follows from one of Lemmas
\ref{lemma:notboth}, \ref{lemma:notbothsmall} or
\ref{lemma:notbothmixed}
    \end{proof}
    
We will also need the following easy lemma 
  
  \begin{lemma} \label{lem:AnexttoB}
 If a configuration carries a label $A^*$ before a saddle move or 
 going through a puncture and a 
 label $B^*$ after then $P_K$ is c-strongly compressible.    
   \end{lemma}
   
   \begin{proof}
       As already discussed the curves before and after the saddle 
       move are distance at most one in the curve complex of $P_K$.
       \end{proof}

\section{Main result} \label{sec:main}
\begin{figure}[tbh]
\centering
\includegraphics[scale=0.45]{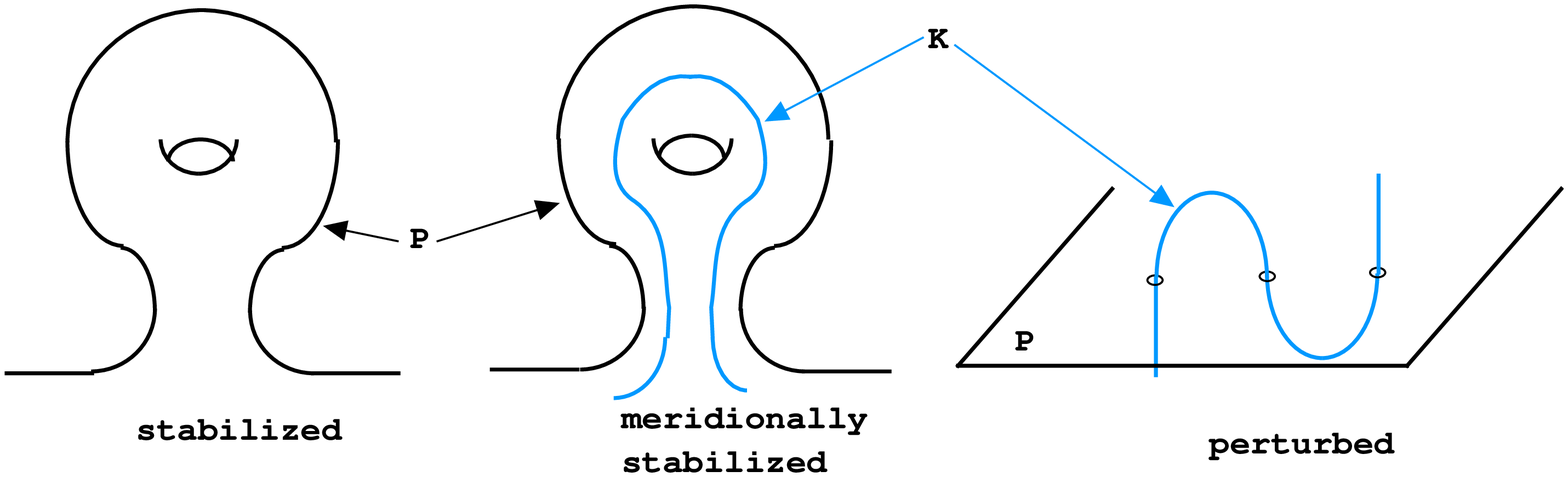}
\caption{} \label{fig:merstab}
\end{figure}

Given a bridge surface for a link $K$ there 
are three ways to create new, more complex, bridge surfaces for the 
link: adding dual one-handles disjoint from the knot (stabilizing), 
adding dual one-handles where one of them has an arc of $K$ as its 
core (meridionally stabilizing), and introducing a pair of a 
canceling minimum and maximum for $K$ (perturbing). These are 
depicted in Figure \ref{fig:merstab}. 

\begin{defin} \label{def:equivalent}
Let $P$ and $Q$ are two bridge surfaces for a knot $K\subset M$. We say 
that $Q$ is equivalent to $P$ if $Q$ is $K$-isotopic to a copy of $P$ 
which may have been stabilized, meridionally stabilized and perturbed.
\end{defin}

There is a fourth way to construct a bridge surface for a knot $K$. 
Suppose $Q$ is a Heegaard splitting for $M$ splitting it into 
handlebodies $X$ and $Y$ and suppose $K$ is isotopic to a subset of 
the spine for $X$. Then by introducing a single minimum, $K$ can be 
placed in bridge position with respect to $Q$. In this case $K$ is 
said to be removable as $Q$ is also a Heegaard surface for $M_K$ 
after an isotopy of $K$. Scharlemann and Tomova discuss all four of these 
operations in detail in \cite{SchTo072}.

Casson and Gordon have demonstrated that if a 3-manifold has a 
Heegaard splitting which is irreducible but strongly compressible 
then the manifold contains an essential surface. In \cite{To072}, Tomova 
extended this result to prove the following

\begin{thm} \label{thm:essentialexists}
Suppose $M$ is a closed orientable irreducible 3-manifold containing 
a link $K$. If $Q$ is a c-strongly compressible bridge surface for 
$K$ then either
\begin{itemize}
    \item $Q$ is stabilized
    \item $Q$ is meridionally stabilized
    \item $Q$ is perturbed
    \item $K$ is removable with respect to $Q$
    \item $M_K$ contains a meridional essential surface $F_K$ such that $2-\chi(F_K) \leq 2-\chi(Q_K)$.
    
\end{itemize}    

\end{thm}

We can now prove the main result of this paper.

\begin{thm} \label{thm:main}
    Suppose $K$ is a nontrivial knot in a closed, irreducible and
    orientable $3$-manifold $M$ and $P$ is a bridge surface for $K$. 
    If $P$ is a sphere assume that $|P \cap K| \geq 6$. If 
    $Q$ is also a bridge surface for $K$ that is not equivalent to 
    $P$, or 
    if $Q$ is a Heegaard surface for $M-\eta (K)$ then
    $d(K,P) \leq 2- \chi(Q-K)$.
    \end{thm}

\begin{proof}
If $Q_K$ is stabilized, meridionally stabilized or perturbed we can perform the
necessary compressions to undo these operations as described in \cite{SchTo072}. Note that these operations increase $\chi(Q_K)$ 
so we may 
assume that $Q_K$ is not stabilized, meridionally stabilized or 
perturbed. If $K$ is removable with respect to $Q$,  we 
may assume that $K$ has been isotoped to lie in the spine of one of 
the handlebodies $M-Q$ so $Q$ is a Heegaard splitting for $M_K$. This operation decreases $|Q \cap K|$ and thus 
also increases $\chi(Q_K)$ .

Suppose first that $Q_K$ is c-strongly
compressible.  If $K$ is not removable with respect to $Q$, by Theorem \ref{thm:essentialexists}, 
there is an
essential surface $F_K$ such that $2-\chi(F_K)<2-\chi(Q_K)$. If $Q$ is 
a Heegaard surface for $M_K$, the existence of such an essential 
surface follows by \cite{CG}. Then the
result follows from Theorem \ref{thm:essbound}.  If $P_K$ is 
c-strongly
compressible, then $d(P,K)\leq 3$ by applying Proposition
\ref{prop:cutimpliescompressing} twice.  Thus we may assume that both $P_K$ and
$Q_K$ are c-weakly incompressible.

The proof now is almost identical to the proof of the 
main result in \cite{ST} so we will only give a brief summary. 

Recall that if $\Sigma_{(A,K)}$ is a spine for the $K$-handlebody $A_K$, 
then $A-\Sigma_{(A,K)} \cong P_K \times I$. Thus if $P$ is a bridge 
surface for $K$, there is a map $H: (P, P\cap K) \times I 
\rightarrow (M,K)$ that is a homeomorphism except over $\Sigma_{(A,K)} 
\cup \Sigma_{(B,K)}$ and near $P \times \bdd I$ the map gives a mapping 
cylinder structure to $\Sigma_{(A,K)} \cup \Sigma_{(B,K)}$. If we 
restrict $H$ to $P_K \times (I, \bdd I) \rightarrow (M, \Sigma_{(A,K)} 
\cup \Sigma_{(B,K)})$ $H$ is called a sweep-out associated to $P$. 

If $Q$ is a Heegaard surface for $M_K$, splitting $M_K$ into a 
compression body and a handlebody, then a similar sweep-out is 
associated to $Q$ between the two spines.  We will denote these spines 
by $\Sss_X$ and $\Sss_Y$.

Consider a square $I \times I$ that describes generic sweep-outs 
of $P_K$ and $Q_K$
from $\Sss_{(A,K)}$ to $\Sss_{(B,K)}$ and from $\Sss_{(X,K)}$ to 
$\Sss_{(Y,K)}$ if $Q$ is a bridge surface for $K$ or from $\Sss_X$ to $\Sss_Y$ if $K$ is removable with 
respect to $Q$.  See Figure 
\ref{fig:graphic}. Each point in the square represents a positioning 
of $P_K$ and $Q_K$.
Inside the square is a graph $\Ggg$, called the {\em graphic} that
represents points at which the intersection is not generic: at each 
point in an edge
in the graphic there is a
single point of tangency between $P_K$ and $Q_K$ or one of the 
surfaces is 
passing through a puncture of the other. At each (valence four) 
vertex of $\Ggg$ there are
two points of tangency or puncture crossings.  By general position 
of, say, the spine
$\Sigma_{(A,K)}$ with the surface $Q_K$ the graphic $\Ggg$ is 
incident to
$\bdd I \times I$ in only a finite number of points (corresponding to
tangencies between $\Sigma_{(A,K)}$ and $Q_K$).  Each such point 
in $\bdd
I \times I$ is incident to at most one edge of $\Ggg$.

Any point in the complement of $\Ggg$ represents a generic 
intersection of $P_K$ and $Q_K$.  Each component of the graphic 
complement will be called a {\em region}; any two points in the same 
region represent isotopic configurations. Label each region 
with labels $\aA,\bB,\xX$ and $\yY$ as described previously where a 
region 
is labeled $\xX$ (resp $\yY$) if any of the labels $x, X, X^c$ (resp 
$y, Y, Y^c$) appear and $\aA$ (resp $\bB$) if the labels $A$ or 
$A^c$ (resp $B$ or $B^*$) appear. See Figure 
\ref{fig:graphic}.
If any region is unlabeled we are 
done by Lemma \ref{lemma:nolabel}.  Also if a region is labeled $\xX$ 
and $\yY$ we are done by one of the Lemmas \ref{lemma:notboth}, 
\ref{lemma:notbothsmall} or
\ref{lemma:notbothmixed}. Finally by Proposition \ref{prop:AandB} no 
region is labeled both $\aA$ and 
$\bB$ so we can assume that each 
region of the square has a unique label.

\begin{figure}[tbh]
\centering
\includegraphics[width=0.6\textwidth]{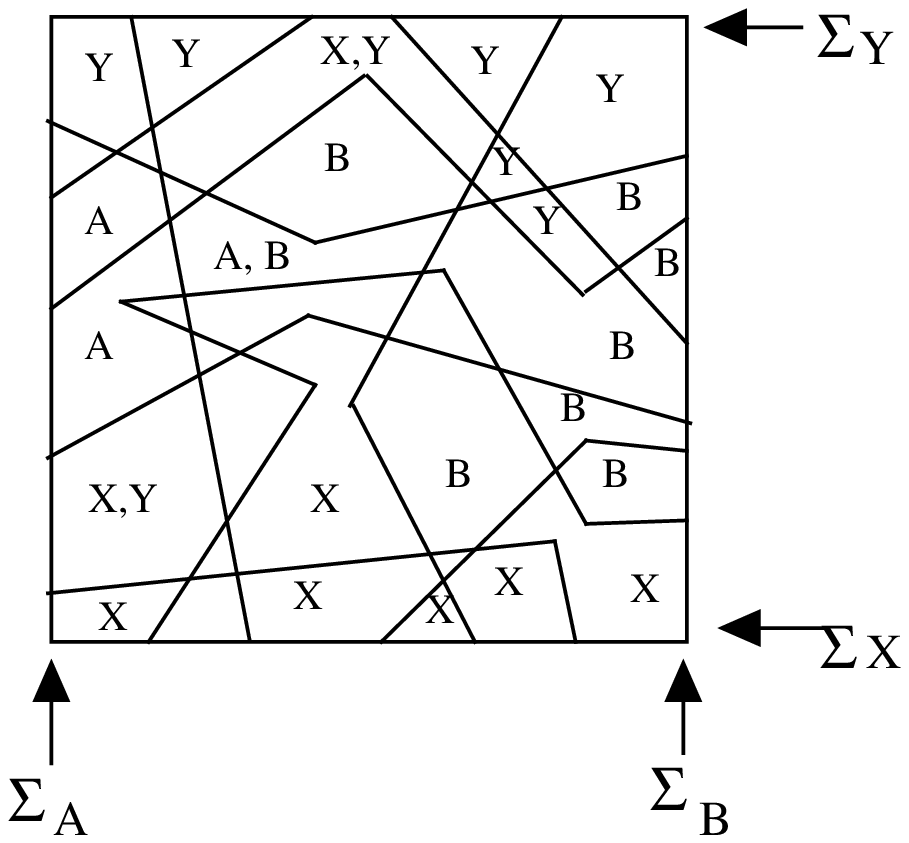}
\caption{} \label{fig:graphic}
\end{figure}

 Let $\Lll$ be the dual complex of $\Ggg$ in $I \times I$; $\Lll$ has one vertex in each face of 
$\Ggg$ and one vertex in each component of $\bdd I \times I - \Ggg$. 
Each edge of $\Lll$ not incident to $\bdd I \times I$ crosses exactly 
one interior edge of $\Ggg$.  Each face of $\Lll$ is a quadrilateral 
and each vertex inherits the label of the corresponding region on 
$\Ggg$.  


Consider the labeling of two adjacent vertices of $\Lll$.
Corollary \ref{corollary:notboth} says that if they are
labeled $\xX$ and $\yY$ we have the desired result and Lemma 
\ref{lem:AnexttoB} says they cannot 
be labeled $\aA$ and $\bB$.  Finally, a discussion identical to the one 
in \cite{ST} about labeling along the edges of $I \times I$ 
shows that no label $\bB$ appears along the $\Sigma_{(A,K)}$ side of $I 
\times 
I$ (the left side in the figure), no label $\aA$ appears along the 
$\Sigma_{(B,K)}$ side (the right side), no label $\yY$ appears along 
the $\Sigma_{(X,K)}$ side ($\Sss_X$ side if $Q$ is a bridge surface 
for $M_K$) (the bottom) and no label $\xX$ appears 
along 
the $\Sigma_{(Y,K)}$ side ($\Sss_Y$ side if $Q$ is a bridge surface 
for $M_K$) (the top).  

We now appeal to the following quadrilateral variant of
Sperner's Lemma proven in the appendix of \cite{ST}:

\begin{lemma} \label{lemma:sperner}  Suppose a square $I \times I$ is 
tiled by quadrilaterals
so that any two that are incident meet either in a corner of each or
in an entire side of each.  Let $\Lll$ denote the graph in $I \times
I$ that is the union of all edges of the quadrilaterals.  Suppose each
vertex of $\Lll$ is labeled $N, E, S,$ or $W$ in such a way that

\begin{itemize}  
    
    \item  no vertex on the East side of $I \times I$ is labeled $W$, 
no vertex
on the West side is labeled $E$, no vertex on the South side is
labeled $N$ and no vertex on the North side is labeled $S$.
    \item no edge in $\Lll$ has ends labeled $E$ and $W$ nor 
    ends labeled $N$ and $S$.
\end{itemize}

Then some quadrilateral contains all four labels

\end{lemma}

In our context the lemma says that there are four regions in 
the graphic incident to the same vertex of $\Ggg$ labeled $\aA, \bB, 
\xX$ and $\yY$. Note then that only two saddle or puncture moves are 
needed to
move from a configuration labeled $\aA$ to one labeled $\bB$.  The
former configuration includes a c-disk for $P_K$ in $A$ and the latter
a c-disk for $P_K$ in $B$.  Note that as $K$ is nontrivial $\chi(Q_K) \leq -2$. Using Proposition 
\ref{prop:cutimpliescompressing} it follows that $d(K,P) \leq 4 
\leq 2 - \chi(Q_K)$, as
long as at least one of the regions labeled $\xX$ and $\yY$ contains
at least one essential curve.

Suppose all curves of $P\cap Q$ in the regions $\xX$ and $\yY$ are
inessential.  Consider the region labeled $\xX$.  Crossing the edge
in the graphic from this region to the region labeled $\aA$ corresponds
to attaching a band $b_A$ with both endpoints on an inessential curve
curve $c \in P\cap Q$ or with endpoint on two distinct curves $c_1$
and $c_2$ where $c_1$ and $c_2$ both bound once punctured disks on
$P_K$.  Note that attaching this band must produce an essential curve
that gives rise to the label $\aA$, call this curve $c_A$.  Similarly
crossing the edge from the region $\xX$ into the region $\bB$
corresponds to attaching a band $b_B$ to give a curve $c_B$.  The two
bands have disjoint interiors and must have at least one endpoint on 
a common curve otherwise $c_A$ and
$c_B$ would be disjoint curves giving rise to labels $\aA$ and $\bB$.  By
our hypothesis attaching both bands simultaneously results in an
inessential curve $c_{AB}$.  We will show that in all cases we can
construct an essential curve $\gamma$ on $P_K$ that is disjoint from
$c_A$ and $c_B$.  After possibly applying Proposition
\ref{prop:cutimpliescompressing}, this implies that $d(K,P)\leq 4$

{\bf Case 1:} Both bands have both of their endpoints on the same curve $c$.

Attaching $b_A$ to $c$ produces two curves that cobound a possibly
once punctured annulus, one of these curves is $c_A$.  We will say
that the band is {\em essential} if $c_A$ is essential on the closed
surface $P$ and {\em inessential} otherwise.  If $b_A$ and $b_B$ are
both essential but $c_{AB}$ is inessential on $P$, then $P$ is a torus
so $P_K$ is a torus with at least two punctures.  In this case $c_A
\cup c_B$ doesn't separate the torus so we can consider the curve
$\gamma$ that bounds a disk on $P$ containing at least two punctures
of $P_K$.

If $b_A$ is essential but $b_B$ isn't, then $c_{AB}$ is parallel to
$c_A$ in $P$ and thus must be essential also so this case cannot
occur.

Finally if both $b_A$ and $b_B$ are inessential on $P$ and $P$ is not
a sphere, then let $\gamma$ be an essential curve on $P$ that is
disjoint from $c_A \cup c_B$.  If $P$ is a sphere, it must have at
least 6 punctures.  Note that $c \cup b_A \cup b_B$ separates $P$ 
into 
4 regions that may contain punctures.  As $P$ has at least 6
punctures, one of these regions contains at least two punctures.  Take
$\gamma$ to be a curve that bounds a disk containing two punctures and
that is disjoint from $c \cup b_A \cup b_B$.

{\bf Case 2:} One band, say $b_A$ has endpoint lying on two different curves
$c_1$ and $c_2$ and the other band, $b_B$ has both endpoints lying on
$c_1$.

If $b_B$ is essential on $P$, then adding both bands simultaneously
results in a curve that is parallel to $c_B$ in $P$ and therefore is
essential contradicting the hypothesis.  If $b_B$ is inessential on
$P$, then $c_1 \cup c_2 \cup b_A \cup b_B$ separates $P$ into 4
regions that may contain punctures. As in the previous case we can
construct an essential curve $\gamma$ on $P_K$ that is disjoint from
$c_A$ and $c_B$ either by taking a curve essential on $P$ or, if $P$ 
is a sphere, by taking 
a curve that lies in one of the four regions and bounds two punctures 
on one side.

{\bf Case 3:} The band $b_A$ has endpoint lying on two different curves
$c_1$ and $c_2$ and $b_B$ has endpoint lying on $c_1$ and $c'_2$,
possibly $c_2=c'_2$.

In this case $c_A$ and $c_B$ are both inessential on $P$ so if $P$ is
not a sphere we can again find a curve $\gamma$ disjoint from both
that is essential on $P$.  If $P$ is a sphere, then $c_1 \cup c_2 \cup
c'_2 \cup b_A \cup b_B$ separates $P$ into 4 regions that may contain
punctures and so we can find a curve $\gamma$ that is essential on
$P_K$ and disjoint from $c_A$ and $c_B$
 as above.

\end{proof}

The curve complex for a 4-times punctured sphere is not connected so a 
bound on the distance of a minimal bridge surface for a 2-bridge knot 
cannot be obtained. However Scharlemann and Tomova have proven the following 
uniqueness result.

\begin{thm} (\cite{SchTo072}, Corollary 4.4) \label{thm:unique}
\label{thm:twobridgeunique}
Suppose $K$ is a knot in $S^3$, 2-bridge with respect to the bridge 
surface $P \cong S^2$, and $K$ is not the unknot. Suppose $Q$ is any 
other bridge surface for $K$. Then either
\begin{itemize}
   \item $Q$ is stabilized
   \item $Q$ is meridionally stabilized
   \item $Q$ is perturbed
   \item $Q$ is properly isotopic to $P$.
\end{itemize}    
\end{thm}
  
\begin{cor}
    Suppose $P$ and $Q$ are two bridge surfaces for a knot $K$ and $K$ 
    is not removable with respect to $Q$. Then either $Q$ is 
    equivalent to $P$ or $d(P)\leq 2-\chi(Q_K)$.
\end{cor}
\begin{proof}
    If $K$ is a two bridge knot with respect to a sphere $P$, then 
    by Theorem \ref{thm:unique}, $Q$ is equivalent to $P$. If $P$ is 
    not a 4 times punctured sphere, the result follows from Theorem 
    \ref{thm:main}.

    \end{proof}

\begin{cor} If $K\subset M^3$ is in bridge position 
   with respect to a Heegaard surface $P$ such 
   that $d(K,P) > 2-\chi(P_K)$ then $K$ has a unique minimal bridge 
position.
   \end{cor}
   
   \begin{proof}
    Suppose $K$ can also be placed in bridge position with respect to 
    a second Heegaard surface $Q$ such that $Q$ is not equivalent to 
    $P$.
    By Theorem \ref{thm:main}, $d(K,P)\leq 2-\chi(Q_K)=2-\chi(P_K)$ contradicting the hypothesis.    
       
    \end{proof}   
    
  \section*{Acknowledgment}
I would like to thank Martin Scharlemann for many helpful
conversations.

     \end{document}